\documentclass[leqno, 12pt]{article}

\usepackage{amsmath}
\usepackage{amsfonts}
\usepackage{amssymb}
\usepackage{amsthm}

\usepackage{booktabs}
\usepackage[table,xcdraw]{xcolor}
\usepackage{multirow}
\usepackage{rotating,tabularx}

\usepackage{scalerel,stackengine}
\stackMath
\newcommand\widecheck[1]{%
\savestack{\tmpbox}{\stretchto{%
  \scaleto{%
    \scalerel*[\widthof{\ensuremath{#1}}]{\kern-.6pt\bigwedge\kern-.6pt}%
    {\rule[-\textheight/2]{1ex}{\textheight}}
  }{\textheight}%
}{0.5ex}}%
\stackon[1pt]{#1}{\scalebox{-1}{\tmpbox}}%
}
\usepackage{scalerel,stackengine}
\stackMath
\newcommand\widebreve[1]{%
\savestack{\tmpbox}{\stretchto{%
  \scaleto{%
    \scalerel*[\widthof{\ensuremath{#1}}]{\kern-.1pt\bigcup\kern-.1pt}  
    {\rule[-\textheight/2]{1ex}{\textheight}}
  }{\textheight}%
}{0.5ex}}%
\stackon[1pt]{#1}{\scalebox{1}{\tmpbox}}%
}
\usepackage{comment}

\def\tod{\buildrel d\over\to}
\def\toL1{\buildrel \mathit{L_1}\over\longrightarrow}
\def\top{\buildrel p\over\to}

\def\towpr1{\buildrel w.pr.1\over\to}
\def\11{\buildrel 1-1\over\longleftrightarrow}
\def\eqd{\buildrel d\over=}

\def\~{\tilde}
\def\stms{\!\setminus\!} 
\def\^{\hat}
\def\u{\breve}
\def\c{\check}

\def\implies{\Rightarrow}

\def\sig{\sigma}
\def\Lam{\Lambda}
\def\Gam{\Gamma}
\def\gam{\gamma}
\def\lam{\lambda}

\def\bet{\beta}

\def\Del{\Delta}
\def\eps{\epsilon}

\def\Ups{\Upsilon}

\def\ome{\omega}
\def\Ome{\Omega}
\def\kap{\kappa}

\def\tr{{\rm tr}}
\def\prtl{\partial}
\def\nul{\emptyset}
\def\E{{\rm E}}
\def\V{{\rm Var}}
\def\C{{\rm Cov}}
\def\R1{{\bf R}^1}
\def\B1{{\cal B}^1}

\def\nid{\noindent}

\def\ts{\textstyle}
\def\t12{\textstyle{1\over2}}

\def\smalltype{\let\rm=\eightrm \let\bf=\eightbf \let\it=\eightit \let\sl=\eightsl
 \baselineskip=9pt \rm}
\font\tenrm=cmr10
\font\tenbf=cmbx10
\font\tenit=cmti10
\font\tensl=cmsl10
\def\medtype{\let\rm=\tenrm \let\bf=\tenbf \let\it=\tenit \let\sl=\tensl
 \baselineskip=12pt \rm}
\font\twelverm=cmr12
\font\twelvebf=cmbx12
\font\twelveit=cmti12 
\font\twelvesl=cmsl12
\def\bigtype{\let\rm=\twelverm \let\bf=\twelvebf \let\it=\twelveit \let\sl=\twelvesl
 \baselineskip=14pt \rm}
 
\begin{document}

\title{How Uniform is the Uniform Distribution on Permutations? }
\author{Michael D. Perlman\footnote{mdperlma@uw.edu. Research supported in part by U.S. Department of Defense Grant H98230-10-C-0263/0000 P0004.}\\Department of Statistics\\
University of Washington
}

\maketitle

\begin{abstract}

\nid For large $q$, does the (discrete) uniform distribution on the set of $q!$ permutations of the vector $(1,2,\dots,q)$ closely approximate the (continuous) uniform distribution  on the $(q-2)$-sphere that contains them? These permutations comprise the vertices of the regular permutohedron, a $(q-1)$-dimensional convex polyhedron. Surprisingly to me, the answer is emphatically no: these permutations are confined to a negligible portion of the sphere, and the regular permutohedron occupies a negligible portion of the ball. However,  $(1,2,\dots,q)$ 
is {\it not} the most favorable configuration for spherical uniformity of permutations. Unlike the permutations of $(1,2,\dots,q)$, the normalized surface area of the largest empty spherical cap among the permutations of the most favorable configuration approaches 0 as $q\to\infty$. Nonetheless, these permutations do not approach spherical uniformity either. 

\end{abstract}

\newpage

{\it This paper is dedicated to the memory of Ingram Olkin, my teacher, mentor, and friend, who introduced so many of us to the joy of majorization.}
\vskip6pt

\nid{\bf 1. Are permutations spherically uniform?}

\nid Column vectors denoted by Roman letters appear in bold type, their components  in plain type; thus ${\bf x}=(x_1,\dots,x_q)'\in\mathbb{R}^q$. For any nonzero ${\bf x}\in\mathbb{R}^q$ ($q\ge2$) let $\Pi({\bf x})$ denote the set of all $q!$ permutations of ${\bf x}$, that is
\begin{equation}\label{Pxdef}
{\Pi({\bf x}})=\{P{\bf x}\mid P\in {\cal P}^q\},
\end{equation}
where ${\cal P}^q$ is the set of all $q\times q$ permutation matrices. In this paper, the following general question is examined:
\vskip4pt

\nid{\it Question 1:  For large $q$, do there exist nonzero vectors ${\bf x}\in\mathbb{R}^q$ such that the (discrete) uniform distribution on $\Pi({\bf x})$ closely approximates  the (continuous) uniform distribution on the $(q-2)$-sphere in which $\Pi({\bf x})$ is contained? Do there exist sequences\footnote{\label{Foot23} Superscripts denote indices, not exponents, unless the contrary is evident.} $\{{\bf x}^q\in\mathbb{R}^q\}$ such that $\Pi({\bf x}^q)$ approaches spherical uniformity as $q\to\infty$? }
\vskip4pt

Because $\Pi({\bf x})$ is invariant under permutations of ${\bf x}$, we may always assume that the components of ${\bf x}$ and ${\bf x}^q$ are ordered, i.e., ${\bf x}, {\bf x}^q\in \mathbb{R}_\le^q$, where
\begin{equation}\label{xorder}
\mathbb{R}_\le^q:=\{{\bf x}\in\mathbb{R}^q\mid x_1\le\cdots\le x_q\}.
\end{equation}
Clearly $\|P {\bf x}\|=\|{\bf x}\|$ for all $P\in{\cal P}^q$, so
\begin{equation}\label{inclusion}
\Pi({\bf x})\subset{\cal S}_{\|{\bf x}\|}^{q-1}\cap {\cal M}_{\bf x}^{q-1},
\end{equation}
where ${\cal S}_\rho^{q-1}$ denotes the $0$-centered $(q-1)$-sphere of radius $\rho$ in $\mathbb{R}^q$ and
\begin{equation}\label{Mtilde}
{\cal M}_{\bf x}^{q-1}:=\{{\bf v}\in\mathbb{R}^q\mid {\bf v}'{\bf e}^q={\bf x}'{\bf e}^q\}
\end{equation}
is the $(q-1)$-dimensional hyperplane containing ${\bf x}$ that is orthogonal to ${\bf e}^q:=(1,\dots,1)'$.
Because ${\cal M}_{\bf x}^{q-1}$ does not contain the origin but we wish to work with 0-centered spheres, we shall translate ${\cal M}_{\bf x}^{q-1}$ to 
\begin{equation}\label{Mo}
\~{\cal M}^{q-1}\equiv\{{\bf v}\in\mathbb{R}^q\mid {\bf v}'{\bf e}^q=0\},
\end{equation}
the $(q\!-\!1)$-dimensional linear subspace parallel to $\!{\cal M}_{\bf x}^{q-1}\!$ and orthogonal to ${\bf e}^q$. 

For this purpose consider the $q\times q$ Helmert orthogonal matrix\begin{align*}
\Gam^{q}&\equiv(\gam_1^q,\gam_ 2^{q},\dots,\gam_q^q)\\
&\equiv\begin{pmatrix}1&1&1&\cdots&1\\
1&-1&1&\cdots&1\\
1&0&-2&\cdots&1\\
\vdots&\vdots&\vdots&\cdots&\vdots\\
 1&0&0&\cdots&-(q-1)\end{pmatrix}
 \begin{pmatrix}\frac{1}{\sqrt{q}}&0&0&\cdots&0\\
 0&\frac{1}{\sqrt{1\cdot2}}&0&\cdots&0\\
  0&0&\frac{1}{\sqrt{2\cdot3}}&\cdots&0\\
    \vdots&\vdots&\vdots&\cdots&\vdots\\
    0&0&0&\cdots&\frac{1}{\sqrt{(q-1)q}}
 \end{pmatrix}\\
& \equiv(\gam_1^{q},\Gam_2^{q}),
\end{align*}
where $\gam_1^{q}\equiv \frac{1}{\sqrt{q}}{\bf e}^q:q\times1$ is the unit vector along the direction of ${\bf e}^q$. By the orthogonality of $\Gam^q$,
\begin{align}
(\Gam_2^q)'\gam_1^q&=0,\label{gam21.1}\\
(\Gam_2^q)'\Gam_2^q&=I_{q-1},\label{gam22.1}\\
\Gam_2^q(\Gam_2^q)'&=I_q-\gam_1^q(\gam_1^q)'=:\Ome_q,\label{gam22.11}
\end{align}
where $I_q$ denotes the $q\times q$ identity matrix. Here $\Ome_q$ is the projection matrix of rank $q-1$ that projects $\mathbb{R}^q$ onto $\~{\cal M}^{q-1}$, so that $\Ome_q {\cal M}_{\bf x}^{q-1}=\~{\cal M}^{q-1}$. 

Let 
${\bf y}$ be the projection of ${\bf x}$ onto $\~{\cal M}^{q-1}$:
\begin{align}
{\bf y}&=\Ome_q {\bf x}={\bf x}-\bar x{\bf e}^q,\label{yOmx}
\end{align}
where 
\begin{equation*}
\bar x\equiv\ts\frac{1}{q}{\bf x}'{\bf e}^q=\ts\frac{1}{q}\sum_{i=1}^qx_i
\end{equation*}
is the average of the $q$ components of ${\bf x}$. Then $\bar y=0$ since ${\bf y}'{\bf e}^q=0$, so
\begin{align}
\ts\|{\bf y}\|^2=\|{\bf x}-\bar x{\bf e}^q\|^2=\sum_{i=1}^q(x_i-\bar x)^2=\sum_{i=1}^q(y_i-\bar y)^2,\label{normy}
\end{align}
which is proportional to their sample variance.
Note that $y_1\le\cdots\le y_q$, so
\begin{equation}\label{yorder}
{\bf y}\in \~{\cal M}_\le^{q-1}:=\~{\cal M}^{q-1}\cap\mathbb{R}_\le^q.
\end{equation}

Because $\Ome_qP=P\Ome_q$ for all $P \in{\cal P}^q$,
\begin{equation}\label{piomega}
\Pi({\bf y})=\Pi(\Ome_q{\bf x})=\Ome_q(\Pi({\bf x})),
\end{equation}
so $\Pi({\bf y})$ is a rigid translation of $\Pi({\bf x})$ and satisfies
\begin{equation}\label{inclusion1}
\Pi({\bf y})\subset{\cal S}_{\|{\bf y}\|}^{q-1}\cap \~{\cal M}^{q-1}=:\~{\cal S}_{\|{\bf y}\|}^{q-2},
\end{equation}
the 0-centered $(q-2)$-sphere of radius $\|{\bf y}\|$ in $\~{\cal M}^{q-1}$. If the (discrete) uniform distribution on $\Pi({\bf y})$ is denoted by $\~{\bf U}_{\bf y}^{q-2}$ and the (continuous) uniform distribution on $\~{\cal S}_{\|{\bf y}\|}^{q-2}$ denoted by $\~{\bf U}_{\|{\bf y}\|}^{q-2}$, then
Question 1 can be restated equivalently as follows:
\vskip4pt


\nid{\it Question 2:  For large $q$, do there exist nonzero vectors ${\bf y}\in \~{\cal M}_\le^{q-1}$ such that $\~{\bf U}_{\bf y}^{q-2}$ closely approximates $\~{\bf U}_{\|{\bf y}\|}^{q-2}$? Do there exist sequences $\{{\bf y}^q\in\~{\cal M}_\le^{q-1}\}$ such that the discrepancy between 
$\~{\bf U}_{{\bf y}^q}^{q-2}$ and $\~{\bf U}_{\|{\bf y}^q\|}^{q-2}$ approaches zero as $q\to\infty$?}
\vskip6pt

\nid{\bf 2. Measures of spherical discrepancy.} 

\nid If we abuse notation by letting $\~{\bf U}_{\bf y}^{q-2}$  and $\~{\bf U}_{\|{\bf y}\|}^{q-2}$ also denote random vectors having these distributions, then the possible existence of the vectors ${\bf y}$ and sequences $\{{\bf y}^q\}$ in Question 2 is supported by the fact that the first and second moments of $\~{\bf U}_{\bf y}^{q-2}$ and $\~{\bf U}_{\|{\bf y}\|}^{q-2}$ coincide:
\begin{align}
\E(\~{\bf U}_{\bf y}^{q-2})&=E(\~{\bf U}_{\|{\bf y}\|}^{q-2})=0,\label{ExpU}\\
\C(\~{\bf U}_{\bf y}^{q-2})&=  \C(\~{\bf U}_{\|{\bf y}\|}^{q-2})=\ts\frac{\|{\bf y}\|^2}{q(q-1)}(qI_q-{\bf e}^q({\bf e}^q)').\label{CovU}
\end{align}
(In fact all odd moments agree since these are 0 by symmetry.) Three measures of the discrepancy between $\~{\bf U}_{\bf y}^{q-2}$ and $\~{\bf U}_{\|{\bf y}\|}^{q-2}$ will be considered.

For nonzero ${\bf w}\in \~{\cal M}^{q-1}$, $-1\le t<1$, and $\rho>0$ define 
\begin{align}
C({\bf w};t)&:=\big\{{\bf v}\in\~{\cal S}_{\|{\bf w}\|}^{q-2}\bigm| {\bf v}'{\bf w}>\|{\bf w}\|^2t\big\},\label{capdef}\\
\~{\cal C}_\rho^{q-2}&:=\big\{C({\bf w};t)\bigm| {\bf w}\in\~{\cal S}_\rho^{q-2},\,-1\le t<1\big\}.\label{capdef1}
\end{align}
Thus $C({\bf w};t)$ is the open spherical cap in $\~{\cal S}_{\|{\bf w}\|}^{q-2}$ of angular half-width $\cos^{-1}(t)$ centered at ${\bf w}$, while $\~{\cal C}_\rho^{q-2}$ is the set of all such spherical caps in $\~{\cal S}_\rho^{q-2}$. 

If ${\bf U}$ is uniformly distributed over the unit $(q-2)$-sphere in $\mathbb{R}^{q-1}$ then for any unit vector ${\bf u}:(q-1)\times1$,
\begin{equation}\label{beta0}
({\bf u}'{\bf U})^2\eqd\mathrm{Beta}\big(\t12,\frac{q-2}{2}\big).
\end{equation}
Thus, if $0\le t<1$ then the normalized $(q-2)$-dimensional surface area of the spherical cap $C({\bf w};t)\subset\~{\cal S}_{\|{\bf w}\|}^{q-2}$ is given by
\begin{align}
\~{\bf U}_{\|{\bf w}\|}^{q-2}(C({\bf w};t))&=\ts\t12\Pr\big[\mathrm{Beta}\big(\frac{1}{2},\frac{q-2}{2}\big)>t^2\big]\label{Betabd}\\
 &=\ts\frac{\Gam(\frac{q-1}{2})}{2\Gam(\frac{1}{2})\Gam(\frac{q-2}{2})}\int_{t^2}^1w^{-\frac{1}{2}}(1-w)^{\frac{q}{2}-2}dw\label{Betabd1}\\
 &=\ts\frac{\Gam(\frac{q-1}{2})}{\Gam(\frac{1}{2})\Gam(\frac{q-2}{2})}\int_t^1(1-v^2)^{\frac{q}{2}-2}dv\label{Betabd2}\\
  &=:\bet^{q-2}(t)\label{betdef}
 \end{align}
 a strictly decreasing smooth function of $t$. 
 
 The following two bounds for $\bet^{q-2}(t)$, $0\le t<1$, will be used. From \eqref{Betabd1},
 \begin{align}
\bet^{q-2}(t)&<\ts\frac{\Gam(\frac{q-1}{2})}{2\Gam(\frac{1}{2})\Gam(\frac{q-2}{2})}\cdot\frac{1}{t}\int_0^{1-t^2}u^{\frac{q}{2}-2}du\nonumber\\
 &=\ts\frac{\Gam(\frac{q-1}{2})}{\Gam(\frac{1}{2})\Gam(\frac{q-2}{2})}\cdot\frac{(1-t^2)^{\frac{q}{2}-1}}{t(q-2)}\nonumber\\
  &\le\ts\sqrt{\frac{q-2}{2\pi}}\cdot\frac{(1-t^2)^{\frac{q}{2}-1}}{t(q-2)}\label{bd100}\\
  &=\ts\frac{(1-t^2)^{\frac{q}{2}-1}}{t\sqrt{2\pi(q-2)}}.\label{bd108}
\end{align}
The inequality used to obtain \eqref{bd100} appears in Wendel [W]. Second, from \eqref{Betabd2} and Wendell's inequality,
 \begin{align}
\ts\t12-\bet^{q-2}(t)
&=\ts\t12\Pr\big[0\le\mathrm{Beta}\big(\frac{1}{2},\frac{q-2}{2}\big)\le t^2\big]\label{ineq5.0}\\
 &=\ts\frac{\Gam(\frac{q-1}{2})}{\Gam(\frac{1}{2})\Gam(\frac{q-2}{2})}\int_0^t(1-v^2)^{\frac{q}{2}-2}dv\nonumber\\
 &\le\ts\frac{\Gam(\frac{q-1}{2})}{\Gam(\frac{1}{2})\Gam(\frac{q-2}{2})}\int_0^te^{-\frac{v^2(q-4)}{2}}dv\nonumber\\
  &=\ts\frac{\Gam(\frac{q-1}{2})\sqrt{2}}{\Gam(\frac{q-2}{2})\sqrt{q-4)}}\cdot\frac{1}{\sqrt{2\pi}}\int_0^{t\sqrt{q-4}}e^{-\frac{z^2}{2}}dz\nonumber\\
  &\le\ts\sqrt{\frac{q-2}{q-4}}\cdot\Big[\Phi\Big(t\sqrt{q-4}\,\Big)-\t12\Big],\label{ineq5}
\end{align}
where $\Phi$ denotes the standard normal cumulative distribution function. 
\vskip4pt

 \nid{\bf Lemma 2.1.} Let $\{t^q\}$ be a sequence in $[0,1)$ and let $0\le\lam\le\infty$. Then
 \begin{equation}\label{doublebd}
\lim_{q\to\infty}\bet^{q-2}(t^q)=1-\Phi(\lam)\iff  \lim_{q\to\infty}t^q\sqrt{q}=\lam.
 \end{equation}

 \nid{\bf Proof.} Let $X_1$ and $X_{q-2}$ denote independent chi-square variates with 1 and $q-2$ degrees of freedom. From \eqref{Betabd} and \eqref{betdef},
 \begin{align*}
 \bet^{q-2}(t^q)&=\ts\t12\Pr\big[\mathrm{Beta}\big(\frac{1}{2},\frac{q-2}{2}\big)> (t^q)^2\big]\\
 &=\ts\t12\Pr\big[\frac{X_1}{X_1+X_{q-2}}>(t^q)^2\big]\\
  &=\ts\t12\Pr\big[\frac{X_1}{X_{q-2}/(q-2)}>\frac{(t^q)^2(q-2)}{1-(t^q)^2}\big].
 \end{align*}
 Thus, because $\frac{X_{q-2}}{q-2}\top1$ by the Law of Large Numbers and $X_1\eqd [N(0,1)]^2$,
  \begin{equation}\label{doublebd2}
\lim_{q\to\infty}\bet^{q-2}(t^q)=1-\Phi(\lam)\iff   \ts\lim\limits_{q\to\infty}\frac{t^q\sqrt{q-2}}{\sqrt{1-(t^q)^2}}=\lam.
 \end{equation}
 It is straightforward to show that
 \begin{equation}\label{tqiff}
 \ts\lim\limits_{q\to\infty}\frac{t^q\sqrt{q-2}}{\sqrt{1-(t^q)^2}}=\lam\iff\lim\limits_{q\to\infty}t^q\sqrt{q}=\lam
 \end{equation}
 (consider the cases $0\le\lam<\infty$ and $\lam=\infty$ separately), hence \eqref{doublebd} holds.\hfill$\square$
\vskip4pt

For nonzero ${\bf y}\in \~{\cal M}_\le^{q-1}$ and any nonempty finite subset $N\subset\~{\cal S}_{\|{\bf y}\|}^{q-2}$, let $\~{\bf U}_N^{q-2}$ denote the (discrete) uniform distribution on $N$; thus $\~{\bf U}_{\bf y}^{q-2}=\~{\bf U}_{\Pi({\bf y})}^{q-2}$. 
\vskip4pt

\nid {\bf Definition 2.2.} 
The {\it normalized spherical cap discrepancy (NSCD)} of $N$
in $\~{\cal S}_{\|{\bf y}\|}^{q-2}$ is defined as\footnote{\label{Foot0} See Leopardi [L1] Def. 2.11.5, Leopardi [L2] \S1, Alishahi and Zamani [AZ] \S1.2.  Unlike [AZ] we divide $|N\cap C|$ by $|N|$ to be able to compare NSCD's of differing dimensions. }
\begin{align}\label{scapdisdef0}
D^{q-2}(N)&:=\sup\ts\Big\{\big|\~{\bf U}_N^{q-2}(C)-\~{\bf U}_{\|{\bf y}\|}^{q-2}(C)\big|\Bigm|C\in\~{\cal C}_{\|{\bf y}\|}^{q-2} \Big\}\\
 &=\sup\ts\Big\{\big|\frac{|N\cap C|}{|N|}-\~{\bf U}_{\|{\bf y}\|}^{q-2}(C)\big|\Bigm| C\in\~{\cal C}_{\|{\bf y}\|}^{q-2}\Big\},\label{scapdisdef}
\end{align}
where $|N\cap C|$ and $|N|$ are the cardinalities of $N\cap C$ and $N$.
The {\it largest empty cap discrepancy (LECD)} of $N$
in $\~{\cal S}_{\|{\bf y}\|}^{q-2}$ is defined as\footnote{\label{Foot00} See [AZ] \S1.2.}
\begin{align}
 L^{q-2}(N)&:=\sup\ts\Big\{\~{\bf U}_{\|{\bf y}\|}^{q-2}(C)\Bigm| C\in\~{\cal C}_{\|{\bf y}\|}^{q-2},\,N\cap C=\nul\Big\}
 \nonumber\\
&=\sup\ts\Big\{\bet^{q-2}(t)\Bigm|\exists{\bf w}\in\~{\cal S}_{\|{\bf y}\|}^{q-2},\,C({\bf w};\,t)\in\~{\cal C}_{\|{\bf y}\|}^{q-2},\,N\cap C({\bf w};\,t)=\nul\Big\}.\ \ \square\label{LECDdef.41}
\end{align}

Obviously
\begin{equation}\label{LleD}
0\le L^{q-2}(N)\le D^{q-2}(N)\le1.
\end{equation}
Note that the suprema in \eqref{scapdisdef0}-\eqref{LECDdef.41} must be maxima, i.e., must be attained. This follows by applying the Blashke Selection Theorem to
\begin{equation*}
\big\{\overline{\mathrm{co}}(C)\bigm| C\in\~{\cal C}_{\|{\bf y}\|}^{q-2}\big\},
\end{equation*}
a collection of closed convex subsets of the closed ball $\~B$ bounded by $\~{\cal S}_{\|{\bf y}\|}^{q-2}$, where $\overline{\mathrm{co}}(C)$ denotes the closed convex hull in $\~B$ of the spherical cap $C$. It follows from this that
\begin{align}
L^{q-2}(N)
&=\bet^{q-2}\big(t(N)\big),\label{LECDdef.2}
\end{align}
where
\begin{equation}
t(N)
=\min\ts\Big\{t\Bigm|\exists\,{\bf w}\in\~{\cal S}_{\|{\bf y}\|}^{q-2},\,C({\bf w};\,t)\in\~{\cal C}_{\|{\bf y}\|}^{q-2},\,N\cap C({\bf w};\,t)=\nul\Big\}.\label{LECDdef.5}
\end{equation}

Define the unit vectors ${\bf z}_k^q\in\~{\cal M}_\le^{q-1}$, $k=1,\dots,q-1$ as follows:
\begin{equation}\label{ak}
{\bf z}_k^q:=\ts\sqrt{\frac{1}{q}}\Big(\underbrace{-\ts\sqrt{\frac{q-k}{k}},\dots,-\sqrt{\frac{q-k}{k}}}_k, \underbrace{\ts\sqrt{\frac{k}{q-k}},\dots,\sqrt{\frac{k}{q-k}}}_{q-k}\Big)'.
\end{equation}
For $1\le k<l\le q-1$, the inner product between ${\bf z}_k^q$ and ${\bf z}_l^q$ is found to be 
\begin{equation}\label{zklinner}
    ({\bf z}_k^q)'{\bf z}_l^q=\sqrt{\ts\frac{k(q-l)}{(q-k)l}}>0.
\end{equation}

\nid{\bf Lemma 2.3.} For nonzero ${\bf y}\in \~{\cal M}_\le^{q-1}$,
\begin{align}
t(\Pi({\bf y})) &=\ts\frac{1}{\|{\bf y}\|}\min\limits_{1\le k\le q-1} {\bf y}'{\bf z}_k^q,
\label{LleD1}\\
L^{q-2}\big(\Pi({\bf y})\big)&\le\ts\frac{1}{2}.\label{LleD2}
\end{align}
\nid{\bf Proof}. For \eqref{LleD1}, it follows from \eqref{capdef} that if ${\bf w}\in\~{\cal S}_{\|{\bf y}\|}^{q-2}$ then
\begin{align}
\Pi({\bf y})\cap C({\bf w};t)=\nul&\iff\max_{P\in{\cal P}^q}(P{\bf y})'{\bf w}\le\|{\bf y}\|^2t.\label{empcapcond}
\end{align}
Thus from \eqref{LECDdef.5} and the Rearrangement Inequality,
\begin{align}
t(\Pi({\bf y}))&=\ts\frac{1}{\|{\bf y}\|^2}\min\limits_{{\bf w}\in\~{\cal S}_{\|{\bf y}\|}^{q-2}}\,\max\limits_{P\in{\cal P}^q}\,(P{\bf y})'{\bf w}\label{empcapcond1}\\
&=\ts\frac{1}{\|{\bf y}\|^2}\min\limits_{{\bf w}\in\~{\cal M}_\le^{q-1},\,\|{\bf w}\|=\|{\bf y}\|}{\bf y}'{\bf w}\label{empcapcond2}
\end{align}
The set $\~{\cal M}_\le^{q-1}$ is a pointed convex simplicial cone\footnote{\label{Last}  The geometric properties of the polyhedral cone $\~{\cal M}_\le^{q-1}$ that we use here stem from its role as a {\it fundamental region} of the finite reflection group (Coxeter group) of all $q\times q$ permutation matrices acting effectively on $\~{\cal M}^{q-1}$. A readable reference is Grove and Benson [GB]; also see Eaton and Perlman [EP].} whose $q-1$ extreme rays are spanned by ${\bf z}_1^q,\dots,{\bf z}_{q-1}^q$, so $\~{\cal M}_\le^{q-1}$ is their nonnegative span. Thus for ${\bf w}\in\~{\cal M}_\le^{q-1}$ with $\|{\bf w}\|=\|{\bf y}\|$,
\begin{equation*}
{\bf w}=\ts\|{\bf y}\|\cdot\frac{\lam_1{\bf z}_1^q+\cdots+\lam_{q-1}{\bf z}_{q-1}^q}{\|\lam_1{\bf z}_1^q+\cdots+\lam_{q-1}{\bf z}_{q-1}^q\|},
\end{equation*}
for some $\lam_1\ge0,\dots,\lam_{q-1}\ge0$ with $\lam_1+\cdots+\lam_{q-1}=1$. Therefore 
\begin{equation*}
{\bf y}'{\bf w}\ge\|{\bf y}\|\min\limits_{1\le k\le q-1} {\bf y}'{\bf z}_k^q,
\end{equation*}
since ${\bf y}\in\~{\cal M}_\le^{q-1}\implies{\bf y}'{\bf z}_k^q\ge0$ by \eqref{zklinner} and $\|\lam_1{\bf z}_1^q+\cdots+\lam_{q-1}{\bf z}_{q-1}^q\|\le1$, hence
\begin{equation}\label{lastineq}
\ts\frac{1}{\|{\bf y}\|^2}\min\limits_{{\bf w}\in\~{\cal M}_\le^{q-1},\,\|{\bf w}\|=\|{\bf y}\|}{\bf y}'{\bf w}\ge \ts\frac{1}{\|{\bf y}\|}\min\limits_{1\le k\le q-1} {\bf y}'{\bf z}_k^q.
\end{equation}
However equality must hold in \eqref{lastineq} because ${\bf w}_k:=\|{\bf y}\|\,{\bf z}_k^q\in\~{\cal M}_\le^{q-1}$ and ${\bf w}_k=\|{\bf y}\|$. This confirms \eqref{LleD1}.

For \eqref{LleD2}, suppose that $L^{q-2}\big(\Pi({\bf y})\big)>\frac{1}{2}$. Then $\Pi({\bf y})$ must be contained in the complement of some closed hemisphere in $\~{\cal S}_{\|{\bf y}\|}^{q-2}$, hence there is some ${\bf v_0}\in\~{\cal S}_{\|{\bf y}\|}^{q-2}$ such that $0>{\bf w}'P{\bf y}$ for all $P\in{\cal P}^q$. Sum over $P$ to obtain $0>{\bf w}'({\bf e}^q({\bf e}^q)'){\bf y}={\bf w}'{\bf e}^q(({\bf e}^q)'{\bf y})=0$, hence a contradiction. \hfill$\square$
\vskip4pt

It is noted in [L1] Lemma 2.11.6 and [L2] \S1 that if $\{N_n\}$ is a sequence of finite sets in $\~{\cal S}_{\|{\bf y}\|}^{q-2}$ ($q$ fixed), then the uniform distribution on $N_n$ converges weakly to $\~{\bf U}_{\|{\bf y}\|}^{q-2}$ as $n\to\infty$ iff $\lim_{n\to\infty} D^{q-2}(N_n)=0$. This motivates the following definition.
\vskip4pt

\nid{\bf Definition 2.4.} A sequence of nonzero vectors $\{{\bf y}^q\in \~{\cal M}_\le^{q-1}\}$ ($q$ varying) is 
{\it asymptotically permutation-uniform (APU)} if 
\begin{equation}\label{ASPUdef}
\lim_{q\to\infty} D^{q-2}(\Pi({{\bf y}^q)})=0;
\end{equation}
it is {\it asymptotically permutation-full (APF)} if 
\begin{equation}\label{ASPFdef}
\hskip68pt\lim_{q\to\infty} L^{q-2}(\Pi({{\bf y}^q)})=0.\hskip130pt\square
\end{equation}
By \eqref{LleD}, APU $\implies$ APF.
\vskip4pt

We also require a definition of asymptotic emptiness for a sequence of nonzero vectors $\{{\bf y}^q\in \~{\cal M}_\le^{q-1}\}$. Because $\Pi({{\bf y}^q)}$ is a finite subset of the sphere $\~{\cal S}_{\|{\bf y}^q\|}^{q-2}$, it always holds that $\~{\cal S}_{\|{\bf y}^q\|}^{q-2}\stms\Pi({{\bf y}^q)}$ is an infinite union of very small empty spherical caps, so a more stringent definition of emptiness is required.
\vskip4pt

\nid{\bf Definition 2.5.} A sequence of nonzero vectors $\{{\bf y}^q\in \~{\cal M}_\le^{q-1}\}$ ($q$ varying) is 
{\it asymptotically permutation-empty (APE)} if $\exists\,\eps>0$ and, for each $q$, $\exists$ a finite collection $\{C_i^q\mid i=1,\dots,n^q\}$ of (possibly overlapping) empty spherical caps in $\~{\cal S}_{\|{\bf y}^q\|}^{q-2}\stms\Pi({{\bf y}^q)}$ such that each $\~{\bf U}_{\|{\bf y}^q\|}^{q-2}(C_i^q)\ge\eps$ and
\begin{equation}\label{limcapunion}
\hskip80pt\ts\lim_{q\to\infty}\~{\bf U}_{\|{\bf y}^q\|}^{q-2}(\cup_{i=1}^{n^q}C_i^q)=1.\hskip120pt\square
\end{equation}

If $\{{\bf y}^q\in \~{\cal M}_\le^{q-1}\}$ is APE then $\Pi({{\bf y}^q)}$ is asymptotically small in the sense that $\Pi({{\bf y}^q)}\subseteq(\cup_{i=1}^{n^q}C_i^q)^c$ with $\~{\bf U}_{\|{\bf y}^q\|}^{q-2}\big((\cup_{i=1}^{n^q}C_i^q)^c\big)\to0$ as $q\to\infty$. That is, $\Pi({{\bf y}^q)}$ occupies only an increasingly negligible portion of the sphere $\~{\cal S}_{\|{\bf y}^q\|}^{q-2}$. Clearly APE $\implies$ not APF $\implies$ not APU.
\vskip2pt

Now modify the definitions of LECD and APF as follows:
\vskip4pt

\nid{\bf Definition 2.6.}  The {\it largest empty cap angular discrepancy (LECAD)} of $N$
in $\~{\cal S}_{\|{\bf y}\|}^{q-2}$ is defined to be
\begin{align}
A^{q-2}(N)
&:=\sup\ts\Big\{\cos^{-1}(t)\Bigm|\exists{\bf w}\in\~{\cal S}_{\|{\bf y}\|}^{q-2},\, C({\bf w};\,t)\in\~{\cal C}_{\|{\bf y}\|}^{q-2},\,N\cap C({\bf w};\,t)=\nul\Big\}\nonumber\\
 &=\cos^{-1}(t(N)),\label{LECDdef.22} 
\end{align}
where $t(N)$ is defined in \eqref{LECDdef.5}.
A sequence of nonzero vectors $\{{\bf y}^q\in \~{\cal M}_\le^{q-1}\}$ ($q$ varying) is {\it asymptotically permutation-dense (APD)} if 
\begin{equation}\label{ASPDdef}
\hskip60pt\lim_{q\to\infty} A^{q-2}(\Pi({{\bf y}^q)})=0.\hskip140pt\square
\end{equation}
Note that \eqref{LECDdef.2} and \eqref{LECDdef.22} yield the relation
\begin{equation}\label{LandA}
L^{q-2}(N)=\bet^{q-2}(\cos(A^{q-2}(N))).
\end{equation}
\vskip2pt

If we set $t^q=t(\Pi({\bf y}^q))$, it follows from \eqref{LECDdef.22} and \eqref{doublebd} with $\lam=\infty$ that 
\begin{equation*}
\{{\bf y}^q\}\  \mathrm{APD}\iff \lim_{q\to\infty}\cos^{-1}(t^q)=0 \iff t^q\to1 \Longrightarrow \bet^{q-2}(t^q)\to0,
\end{equation*}
hence APD $\implies$ APF. However the converse need not hold: it will be shown in \S4 that the sequence $\{\^{\bf y}^q\}$ of {\it maximal configurations} defined in \eqref{yhatq} is APF but not APD.
\vskip4pt

\nid{\bf Remark 2.7.} Consider a sequence of spherical caps $C({\bf w}^q;t^q)\subseteq\~{\cal M}^{q-1}$ such that $t^q\to0$ while $t^q\sqrt{q}\to\infty$. Then $\cos^{-1}(t^q)\to\frac{\pi}{2}$, while $\bet^{q-2}(t^q)\to0$ by \eqref{doublebd} with $\lam=\infty$, that is, the spherical caps approach hemispheres in terms of their angular measure but their surface areas approach 0. An example can be seen in \S4 by taking $C({\bf w}^q;t^q)$ to be the largest empty spherical cap for the set $\Pi(\^{\bf y}^q)$, see \eqref{Lhatlimit} and \eqref{Ahatlimit4}.\hfill$\square$

\vskip4pt
Question 2 now can be refined further as follows:

\vskip4pt

\nid{\it Question 3:  For which ${\bf y}\in \~{\cal M}_\le^{q-1}$, if any, are $D^{q-2}(\Pi({\bf y}))$, $L^{q-2}(\Pi({\bf y}))$, and/or $A^{q-2}(\Pi({\bf y}))$ small? Which sequences $\{{\bf y}^q\in \~{\cal M}_\le^{q-1}\}$, if any, are APU? APF? APD? APE?}
\vskip4pt

Some answers to these questions will be derived in \S3-\S5 and summarized in \S6; for example, no APD sequence exists (Proposition 6.1). Some results about the volumes of the corresponding permutohedra with vertices $\Pi({\bf y}^q)$ are presented in \S7. Several open questions are stated in \S6-8.
\vskip4pt

\nid {\bf Example 2.8.} Despite the agreement of the first and second moments of $\~{\bf U}_{\bf y}^{q-2}$ and $\~{\bf U}_{\|{\bf y}\|}^{q-2}$ (cf. \eqref{ExpU}, \eqref{CovU}), $L^{q-2}(\Pi({\bf y}))$ need not be small. For example, take ${\bf y}={\bf f}_q^q$ where, for $i=1,\dots,q$, ${\bf f}_i^q\in\~{\cal M}_\le^{q-1}$ is the unit column vector
\begin{equation}\label{fi}
{\bf f}_i^q=\ts\frac{1}{\sqrt{q(q-1)}}(\underbrace{-1,\dots,-1}_{i-1},\,q-1,\,\underbrace{-1\dots,-1}_{q-i})'.
\end{equation}
Here $\Pi({{\bf f}_q^q})=\{{\bf f}_1^q,\dots,{\bf f}_q^q\}$, so $|\Pi({{\bf f}_q^q})|=q$ not $q!$. From \eqref{LleD1}, \eqref{LECDdef.2}, and \eqref{doublebd} with $t^q=\frac{1}{q-1}$ and $\lam=0$,
\begin{align}
t(\Pi({{\bf f}_q^q}))&=\ts\frac{1}{q-1},\label{fqq1}\\
L^{q-2}(\Pi({{\bf f}_q^q}))&=\ts\bet^{q-2}\big(\frac{1}{q-1}\big)\to\t12\label{fqq2}
\end{align}
as $q\to\infty$. Thus the sequence $\{{\bf f}_q^q\}$ is not APF, hence not APU. \hfill $\square$
\vskip4pt

\nid{\bf Remark 2.9.} For later use, we note that for $i=1,\dots,q$,
\begin{equation}\label{fidef}
{\bf f}_i^q=\ts\frac{1}{\sqrt{q(q-1)}}(q{\bf e}_i^q-{\bf e}^q)=\sqrt{\frac{q}{q-1}}\Ome_q {\bf e}_i^q,
\end{equation}
where ${\bf e}_i^q\equiv(0,\dots,0,1,0,\dots,0)'$ denotes the $i$th coordinate vector in $\mathbb{R}^q$ and
 $\Ome_q {\bf e}_i^q$ is the projection of ${\bf e}_i^q$ onto $\~{\cal M}^{q-1}$. Thus ${\bf f}_1^q,\dots,{\bf f}_q^q$ form the vertices of a standard simplex in $\~{\cal M}^{q-1}$: an equilateral triangle when $q=3$, a regular tetrahedron when $q=4$, etc.\hfill $\square$
\vskip4pt

For any nonzero ${\bf y}^q\in \~{\cal M}_\le^{q-1}$, 
$\frac{\sqrt{q-1}}{\|{\bf y}^q\|}\~{\bf U}_{\|{\bf y}^q\|}^{q-2}$ is uniformly distributed on the sphere of radius $\sqrt{q-1}$ in $\~{\cal M}^{q-1}$. It is well known (e.g. Eaton [E] Proposition 7.5), and also follows from \eqref{Betabd}-\eqref{betdef} and Lemma 2.1, that the marginal distributions from this uniform distribution converge to the standard normal distribution $N(0,1)$ as $q\to\infty$. More precisely, for any sequence of unit vectors $\{{\bf u}^q\}$ in $\~{\cal M}^{q-1}$,
\begin{align*}
\ts({\bf u}^q)'\Big(\frac{\sqrt{q-1}}{\|{\bf y}^q\|}\~{\bf U}_{\|{\bf y}^q\|}^{q-2}\Big)=\frac{\sqrt{q-1}}{\|{\bf y}^q\|}({\bf u}^q)'\~{\bf U}_{\|{\bf y}^q\|}^{q-2}
&\tod N(0,1).
\end{align*}
as $q\to\infty$. If we take ${\bf u}^q={\bf f}_i^q=\sqrt{\frac{q}{q-1}}\Ome_q{\bf e}_i^q$ (see \eqref{fidef}) for any fixed $i$, where ${\bf e}_i^q\equiv(0,\dots,0,1,0,\dots,0)'$ is the $i$th coordinate vector in $\mathbb{R}^q$, then
\begin{equation}\label{coordNormal}
\ts\frac{\sqrt{q-1}}{\|{\bf y}^q\|}({\bf u}^q)'\~{\bf U}_{\|{\bf y}^q\|}^{q-2}=\frac{\sqrt{q}}{\|{\bf y}^q\|}({\bf e}_i^q)'\~{\bf U}_{\|{\bf y}^q\|}^{q-2}\equiv\frac{\sqrt{q}}{\|{\bf y}^q\|}(\~{\bf U}_{\|{\bf y}^q\|}^{q-2})_i\tod N(0,1)
\end{equation}
as $q\to\infty$, where $(\~{\bf U}_{\|{\bf y}^q\|}^{q-2})_i$ denotes the $i$th component of $\~{\bf U}_{\|{\bf y}^q\|}^{q-2}$. 
\vskip4pt

\nid{\bf Proposition 2.10.} A necessary condition that a sequence $\{{\bf y}^q\in \~{\cal M}_\le^{q-1}\}$ of nonzero vectors be APU is that for each fixed $i\ge1$,
\begin{equation}\label{necCondAPU}
\ts\frac{\sqrt{q}}{\|{\bf y}^q\|}(\~{\bf U}_{{\bf y}^q}^{q-2})_i\tod N(0,1)
\end{equation}
as $q\to\infty$, where $(\~{\bf U}_{{\bf y}^q}^{q-2})_i$ denotes the $i$th component of $\~{\bf U}_{{\bf y}^q}^{q-2}$.
\vskip2pt
\nid{\bf Proof.} From \eqref{scapdisdef0} and \eqref{capdef}-\eqref{capdef1},
\begin{align*}
&D^{q-2}(\Pi({\bf y}^q))\\
&=\sup\ts\big\{\big|\~{\bf U}_{{\bf y}^q}^{q-2}(C)-\~{\bf U}_{\|{\bf y}^q\|}^{q-2}(C)\big|\bigm|C\in\~{\cal C}_{\|{\bf y}\|}^{q-2} \big\}\\
&\ge\sup_{-1\le t<1}\ts\big|\~{\bf U}_{{\bf y}^q}^{q-2}(C(\|{\bf y}^q\|{\bf f}_i^q;t))-\~{\bf U}_{\|{\bf y}^q\|}^{q-2}(C(\|{\bf y}^q\|{\bf f}_i^q;t))\big|\\
&=\sup_{-1\le t<1}\big|\Pr\!\big[({\bf f}_i^q)'\~{\bf U}_{{\bf y}^q}^{q-2}>\|{\bf y}^q\|\,t\big]-\Pr\!\big[({\bf f}_i^q)'\~{\bf U}_{\|{\bf y}^q\|}^{q-2}>\|{\bf y}^q\|\,t\big]\big|\\
&=\ts\sup\limits_{-1\le t<1}\big|\Pr\!\big[\sqrt{\frac{q}{q-1}}({\bf e}_i^q)'\~{\bf U}_{{\bf y}^q}^{q-2}>\|{\bf y}^q\|\,t\big]-\Pr\!\big[\sqrt{\frac{q}{q-1}}({\bf e}_i^q)'\~{\bf U}_{\|{\bf y}^q\|}^{q-2}>\|{\bf y}^q\|\,t\big]\big|\\
&=\ts\sup\limits_{-1\le t<1}\big|\Pr\!\big[\frac{\sqrt{q}}{\|{\bf y}^q\|}(\~{\bf U}_{{\bf y}^q}^{q-2})_i>\sqrt{q-1}\,t\big]-\Pr\!\big[\frac{\sqrt{q}}{\|{\bf y}^q\|}(\~{\bf U}_{\|{\bf y}^q\|}^{q-2})_i>\sqrt{q-1}\,t\big]\big|.
\end{align*}
Because $D^{q-2}(\Pi({\bf y}^q))\to0$ if $\{{\bf y}^q\}$ is APU, this and \eqref{coordNormal} yield \eqref{necCondAPU}.  
\hfill$\square$
 \vskip6pt
 

\nid{\bf 3. The regular configurations $\bar{\bf x}^q$ and $\bar{\bf y}^q$ are not spherically uniform.} 

\nid It is seen from \eqref{LleD2} and \eqref{fqq2} that $\{{\bf f}_q^q\}$ fails to be APF (hence fails to be APU and APD) to the greatest possible extent. Clearly this is due to the fact that the components of ${\bf f}_q^q$ comprise only two distinct values $-1$ and $q-1$. This suggests that the APU, APF, and APD properties are more likely to hold for vectors ${\bf y}^q\equiv\Ome_q{\bf x}^q\in\~{\cal M}_\le^{q-1}$ whose components are distinct, so that $|\Pi({\bf y}^q)|$, equivalently $|\Pi({\bf x}^q)|$, attains its maximum value $q!$. 

At this point, it seems reasonable to conjecture that the APU, APD, and APF properties are most likely to hold for vectors whose components are evenly spaced, that is, for the vectors
\begin{align}
\bar{\bf x}^q&=(1,2,\dots,q)',\label{xqdef}\\
\bar{\bf y}^q=\Ome_q \bar{\bf x}^q&=\ts(-\frac{q-1}{2},-\frac{q-3}{2},\dots,\frac{q-3}{2},\frac{q-1}{2})',\label{yqdef}
\end{align}
We call $\bar{\bf x}^q$ and $\bar{\bf y}^q$ the {\it regular configurations} in $\mathbb{R}_\le^q$ and $\~{\cal M}_\le^{q-1}$ respectively. 

This conjecture is supported by the case $q=2$ with $\bar{\bf y}_2=(-\frac{1}{2},\,\frac{1}{2})'$, where the two permutations $(-\frac{1}{2},\,\frac{1}{2})'$ and $(\frac{1}{2},\,-\frac{1}{2})'$ trivially are uniformly distributed on $\~{\cal S}_{ \|\bar{\bf y}^2\|}^0$, and by the case $q=3$ with $\bar{\bf y}^3=(1,2,3)'$, where the 3!=6 permutations of $\bar{\bf y}^3$ comprise the vertices of a regular hexagon, the most uniform among all configurations of 6 points on the circle $\~{\cal S}_{\|\bar{\bf y}^3\|}^1$. When $q=4$, however, the 4!=24 permutations of ${\bf y}^{(4)}\equiv(-\frac{3}{2},-\frac{1}{2},\frac{1}{2},\frac{3}{2})$ comprise the vertices of the regular permutohedron $\~{\mathfrak{R}}^4$ (see \S7), a truncated octahedron whose 14 faces consist of 8 regular hexagons and 6 squares, hence is not a regular solid.

In this section we present two arguments that show this asymptotic spherical uniformity conjecture is invalid for the regular configurations. The first argument (Propositions 3.1 and 3.2) examines the APF and APE properties for $\{\bar{\bf x}^q\}$ and $\{\bar{\bf y}^q\}$, the second argument (Proposition 3.4) compares the univariate marginal distributions of $\~{\bf U}_{\bar{\bf y}^q}^{q-2}$  and $\~{\bf U}_{\|\bar{\bf y}^q\|}^{q-2}$.  A third comparison of $\~{\bf U}_{\bar{\bf y}^q}^{q-2}$ and $\~{\bf U}_{\|\bar{\bf y}^q\|}^{q-2}$ will be presented in \S7.

\vskip4pt

\nid{\bf Proposition 3.1.} The sequences of regular configurations $\{\bar{\bf x}^q\in\mathbb{R}_\le^q\}$ and $\{\bar{\bf y}^q\in\~{\cal M}_\le^{q-1}\}$ are not APF, hence not APU and not APD.
\vskip2pt

\nid{\bf Proof.} It suffices to consider $\{\bar{\bf y}^q\}$. Beginning with the relations
\begin{align}
\|\bar{\bf y}^q\|&=\ts\sqrt{\frac{1}{12}q(q^2-1)},\label{normyq}\\
(\bar{\bf y}^q)'{\bf z}_k^q&=(\bar{\bf x}^q)'{\bf z}_k^q=\ts\frac{1}{2}\sqrt{qk(q-k)},
\end{align}
it follows from \eqref{LECDdef.2} and \eqref{LleD1} that the LECD of $\Pi(\bar{\bf y}^q)$ is given by
\begin{align}
L^{q-2}\big(\Pi(\bar{\bf y}^q)\big)&=\ts\bet^{q-2}\Big(\frac{1}{\|\bar{\bf y}^q\|}\min\limits_{1\le k\le q-1} (\bar{\bf y}^q)'{\bf z}_k^q\Big)\label{LleD111}\\
&=\ts\bet^{q-2}\Big(\sqrt{\frac{3}{q+1}}\,\Big),\label{LleD211}
\end{align}
where the minimum is attained for $k=1$ and $k=q-1$. From Lemma 2.1 with $\lam=\sqrt{3}$,
\begin{align}
\ts\lim\limits_{q\to\infty}\bet^{q-2}\Big(\sqrt{\frac{3}{q+1}}\,\Big)
&=1-\Phi(\sqrt{3})\approx.0416>0,\label{ineq6}
\end{align}
so $\{\bar{\bf y}^q\}$ is not APF.\hfill$\square$
\vskip4pt

In fact, $\{\bar{\bf x}^q\}$ and $\{\bar{\bf y}^q\}$ fail asymptotic uniformity in a stronger sense:
\vskip4pt

\nid{\bf Proposition 3.2.} The regular configurations $\{\bar{\bf x}^q\}$
and $\{\bar{\bf y}^q\}$
are APE.
\vskip2pt

\nid{\bf Proof.} Again it suffices to consider $\{\bar{\bf y}^q\}$. Define $\bar{\bf z}_k^q=\|\bar{\bf y}^q\|{\bf z}_k^q$, where ${\bf z}_k^q$ is the unit vector in \eqref{ak}. Because the minimum in \eqref{LleD111} is attained for $k=1$ and $q-1$, i.e., for ${\bf z}_1^q$ and ${\bf z}_{q-1}^q$, both $C\big(\bar{\bf z}_1^q;\sqrt{\frac{3}{q+1}}\big)$ and $C\big(\bar{\bf z}_{q-1}^q;\sqrt{\frac{3}{q+1}}\big)$ are (overlapping) largest empty spherical caps for $\Pi(\bar{\bf y}^q)$ in $\~{\cal S}_{\|\bar{\bf y}^q\|}^{q-2}$. (Note that ${\bf z}_1^q=-{\bf f}_1^q$ and ${\bf z}_{q-1}^q={\bf f}_q^q$.) Because $P\Pi(\bar{\bf y}^q)=\Pi(\bar{\bf y}^q)$ for all $P\in{\cal P}^q$, $C\big(P\bar{\bf z}_1^q;\sqrt{\frac{3}{q+1}}\big)$ and $C\big(P\bar{\bf z}_{q-1}^q;\sqrt{\frac{3}{q+1}}\big)$ also are (overlapping) largest empty spherical caps for $\Pi(\bar{\bf y}^q)$; there are $2q!$ such caps, all congruent. However
\begin{align*}
\{P\bar{\bf z}_1^q\mid P\in{\cal P}^q\}&=\{-\bar{\bf f}_1^q,\dots,-\bar{\bf f}_q^q\},\\
\{P\bar{\bf z}_{q-1}^q\mid P\in{\cal P}^q\}&=\{\bar{\bf f}_1^q,\dots,\bar{\bf f}_q^q\},
\end{align*}
where $\bar{\bf f}_i^q=\|\bar{\bf y}^q\|{\bf f}_i^q$, so these $2q!$ empty caps reduce to $2q$, namely
\begin{equation*}
\ts\Big\{C\big(-\bar{\bf f}_i^q;\sqrt{\frac{3}{q+1}}\big)\Bigm| i=1,\dots q\Big\}\cup \Big\{C\big(\bar{\bf f}_i^q;\sqrt{\frac{3}{q+1}}\big)\Bigm| i=1,\dots q\Big\}.
\end{equation*}
By \eqref{LleD111}-\eqref{ineq6}, each of these congruent empty caps remains nonnegligible as $q\to\infty$, so $\{\bar{\bf y}^q\}$ is APE if
\begin{equation*}
\lim_{q\to\infty}\~{\bf U}_{\|\bar{\bf y}^q\|}^{q-2}(\Ups^q)=1,
\end{equation*}
where
\begin{align*}
\Ups^q&=\ts\bigcup\nolimits_{i=1}^q \Big[C\big(-\bar{\bf f}_i^q;\sqrt{\frac{3}{q+1}}\big)\cup C\big(\bar{\bf f}_i^q;\sqrt{\frac{3}{q+1}}\big)\Big].
\end{align*}
Therefore, because
\begin{equation*}
\ts\~{\cal S}_{\|\bar{\bf y}^q\|}^{q-2}\bigcap\big((\Ups_q)^c\big)=\~{\cal S}_{\|\bar{\bf y}^q\|}^{q-2}\bigcap\big(\cap_{i=1}^q S_i^q\big),
\end{equation*}
where $S_i^q$ is the closed symmetric slab
\begin{align*}
S_i^q
 &=\ts\Big\{{\bf v}\in\~{\cal M}^{q-1}\bigm| |{\bf v}'\bar{\bf f}_i^q|\le\|\bar{\bf y}^q\|^2\sqrt{\frac{3}{q+1}}\Big\},
\end{align*}
to show that $\{\bar{\bf y}^q\}$ is APE it suffices to show that
\begin{equation}\label{limUslab}
\lim_{q\to\infty}\~{\bf U}_{\|\bar{\bf y}^q\|}^{q-2}\big((\Ups_q)^c\big)\equiv\lim_{q\to\infty}\~{\bf U}_{\|\bar{\bf y}^q\|}^{q-2}\big(\cap_{i=1}^q S_i^q\big)=0.
\end{equation}

If $S_1^q,\dots,S_q^q$ were mutually geometrically orthogonal, i.e., if ${\bf f}_1^q,\dots,{\bf f}_q^q$ were orthonormal, then the $S_i^q$ would be subindependent under $\~{\bf U}_{\|\bar{\bf y}^q\|}^{q-2}$ (cf. Ball and Perissinaki [BP]), that is,
\begin{equation*}
\~{\bf U}_{\|\bar{\bf y}^q\|}^{q-2}(\cap_{i=1}^q S_i^q)\le\prod\nolimits_{i=1}^q\~{\bf U}_{\|\bar{\bf y}^q\|}^{q-2}(S_i^q),
\end{equation*}
which would readily yield \eqref{limUslab}. However, $({\bf f}_i^q)'{\bf f}_j^q=-\frac{1}{q-1}\ne0$ if $i\ne j$ so this approach fails.\footnote{\label{Foot60} In fact, Theorem 2.1 of Das Gupta {\it et al.} [DEOPSS] suggests that $S_1^q,\dots,S_q^q$ may be superdependent under $\~{\bf U}_{\|\bar{\bf y}^q\|}^{q-2}$.} Instead we can apply the cruder one-sided bound
\begin{align}\label{Honesidedbd}
\~{\bf U}_{\|\bar{\bf y}^q\|}^{q-2}(\cap_{i=1}^q S_i^q)
&\le\~{\bf U}_{\|\bar{\bf y}^q\|}^{q-2}(\cap_{i=2}^q H_i^q),
\end{align}
 where $H_i^q$ is the halfspace
\begin{align*}
H_i^q&:=\ts\big\{{\bf v}\in\~{\cal M}^{q-1}\bigm| {\bf v}'\bar{\bf f}_i^q\le\|\bar{\bf y}^q\|^2\sqrt{\frac{3}{q+1}}\big\}.
\end{align*}
Again $H_2^q,\dots,H_q^q$ are not mutually geometrically orthogonal, but now this works in our favor: because $({\bf f}_i^q)'{\bf f}_j^q<0$ if $i\ne j$, the extension of Slepian's inequality to spherically symmetric density functions ([DEOPSS], Lemma 5.1) and a standard approximation argument yields
\begin{equation}\label{aaa}
\~{\bf U}_{\|\bar{\bf y}^q\|}^{q-2}(\cap_{i=2}^q H_i^q)\le\~{\bf U}_{\|\bar{\bf y}^q\|}^{q-2}(\cap_{i=2}^qK_i^q ),
\end{equation}
where $K_i^q$ is the halfspace
\begin{equation*}
\ts K_i^q:=\big\{{\bf v}\in\~{\cal M}^{q-1}\bigm| {\bf v}'\gam_i^q\le \|\bar{\bf y}^q\|\sqrt{\frac{3}{q+1}} \big\}
\end{equation*}
and $\gam_2^q,\dots,\gam_q^q$ are the last $q-1$ columns of the Helmert matrix $\Gam^q$ in \S1, which form an orthonormal basis in $\~{\cal M}^{q-1}$ so $(\gam_i^q)'\gam_j^q=0$.  Now Proposition A.1 in the Appendix and the orthogonal invariance of $\~{\bf U}_{\|\bar{\bf y}^q\|}^{q-2}$ imply that 
\begin{align}
\~{\bf U}_{\|\bar{\bf y}^q\|}^{q-2}(\cap_{i=2}^qK_i^q)&\le\ts\prod_{i=2}^q\~{\bf U}_{\|\bar{\bf y}^q\|}^{q-2}(K_i^q)\\
 &=\ts \big[\~{\bf U}_{\|\bar{\bf y}^q\|}^{q-2}(K_i^q)\big]^{q-1}\\
  &=\ts\Big[1-\bet^{q-2}\Big(\sqrt{\frac{3}{q+1}}\Big)\Big]^{q-1}.
\end{align}
Therefore by \eqref{ineq6},
\begin{equation}
\ts\limsup\limits_{q\to\infty} \big[\~{\bf U}_{\|\bar{\bf y}^q\|}^{q-2}(\cap_{i=2}^qK_i^q)\big]^{\frac{1}{q-1}}
\le\Phi(\sqrt{3})\approx.9584, 
\end{equation}
hence by \eqref{limUslab}-\eqref{aaa},
\begin{equation}\label{Ugeometric}
\~{\bf U}_{\|\bar{\bf y}^q\|}^{q-2}((\Ups^q)^c)\le\~{\bf U}_{\|\bar{\bf y}^q\|}^{q-2}(\cap_{i=2}^qK_i^q)\le(.96)^{q-1}
\end{equation}
for sufficiently large $q$. Thus \eqref{limUslab} holds, in fact $\~{\bf U}_{\|\bar{\bf y}^q\|}^{q-2}((\Ups^q)^c)\to0$ at a geometric rate, hence $\{\bar{\bf y}^q\}$ is APE as asserted.\hfill$\square$
\vskip4pt

\nid{\bf Remark 3.3.} The above result can be framed in terms of statistical hypothesis testing. Based on one random observation ${\bf Y}\equiv(Y_1,\dots,Y_q)'\in\~{\cal S}_{\|\bar{\bf y}^q\|}^{q-2}$, suppose that it is wished to test the spherical-uniformity hypothesis $H_0$ that ${\bf Y}\eqd\~{\bf U}_{\|\bar{\bf y}^q\|}^{q-2}$ against the permutation-uniformity alternative $H_1$ that ${\bf Y}\eqd\~{\bf U}_{\bar{\bf y}^q}^{q-2}$. Consider the test that rejects $H_0$ in favor of $H_1$ iff ${\bf Y}\in(\Ups^q)^c$, that is, iff 
\begin{equation*}
\ts\max_{1\le i\le q}|Y_i-\bar Y|\le\frac{q-1}{2},
\end{equation*}
where $\bar Y=\frac{1}{q}\sum_{i=1}^q Y_i$. The size of this test is $\~{\bf U}_{\|\bar{\bf y}^q\|}^{q-2}((\Ups^q)^c)$, which by \eqref{Ugeometric} rapidly approaches 0 as $q\to\infty$, while its power = 1 for every $q$ because $\Pi(\bar{\bf y}^q)\subset(\Ups^q)^c$.\hfill$\square$
\vskip4pt

A second argument for the invalidity of the spherical uniformity conjecture for the regular configuration $\{\bar{\bf y}^q\}$ (and $\{\bar{\bf x}^q\}$) stems from Proposition 2.10 and the following fact:
\vskip4pt

\nid{\bf Proposition 3.4.} For each fixed $i\ge1$, as $q\to\infty$,
\begin{align}
\ts\sqrt{\frac{12}{q^2-1}}(\~{\bf U}_{\bar{\bf y}^q}^{q-2})_i&\tod\mathrm{Uniform}\big(-\sqrt{3},\sqrt{3}\,\big)\label{Uniformlim}
\end{align}
as $q\to\infty$. Thus $\{\bar{\bf y}^q\}$ does not satisfy \eqref{necCondAPU}, hence is not APU.
\vskip2pt

\nid{\bf Proof.} By \eqref{yqdef}, for each $i=1,\dots,q$, $(\~{\bf U}_{\bar{\bf y}^q}^{q-2})_i$ is uniformly distributed over the range
\begin{equation}\label{rangeU}
\ts-\frac{q-1}{2},-\frac{q-3}{2},\dots,\frac{q-3}{2},\frac{q-1}{2},
\end{equation}
so its moment generating function (mgf) is
\begin{equation*}
\frac{e^{tq/2}-e^{-tq/2}}{q(e^{t/2}-e^{-t/2})}.
\end{equation*}
 (Thus the distribution of $(\~{\bf U}_{\bar{\bf y}^q}^{q-2})_i$ is the same for each $i$.) Therefore the mgf of $\ts\sqrt{\frac{12}{q^2-1}}(\~{\bf U}_{\bar{\bf y}^q}^{q-2})_i$ is
\begin{equation*}
\frac{e^{t\sqrt{\frac{3q^2}{q^2-1}}}-e^{-t\sqrt{\frac{3q^2}{q^2-1}}}}{q\Big(e^{t\sqrt{\frac{3}{q^2-1}}}-e^{-t\sqrt{\frac{3}{q^2-1}}}\Big)},
\end{equation*}
which converges to $\frac{\sinh(t\sqrt{3})}{t\sqrt{3}}$ as $q\to\infty$, the mgf of $\mathrm{Uniform}\big[-\sqrt{3},\sqrt{3}\,\big]$.
\hfill$\square$
\vskip6pt

\nid{\bf 4. The most favorable configuration for spherical uniformity.} 

\nid It was shown in Proposition 3.1 that the regular configurations $\bar{\bf x}^q$ and $\bar{\bf y}^q$ are not APF, hence not APU or APD, although the components of $\bar{\bf x}^q$ and $\bar{\bf y}^q$ are exactly evenly spaced.  Is there is a more favorable configuration for spherical uniformity of permutations? We show now that the answer is yes.

Continuing the discussion in \S2-3, we wish to find a nonzero vector ${\bf y}$ in $\~{\cal S}_{\|\bar{\bf y}^q\|}^{q-2}\cap\mathbb{R}_\le^q$ that minimizes the LECD $L^{q-2}(\Pi({\bf y}))$ in $\~{\cal S}_{\|\bar{\bf y}^q\|}^{q-2}$;
equivalently, that minimizes the LECAD $A^{q-2}(\Pi({\bf y}))$. From \eqref{LECDdef.2}, \eqref{LleD1}, and \eqref{LandA},
\begin{align}
L^{q-2}(\Pi({\bf y}))&=\bet^{q-2}(t(\Pi({\bf y})))=\ts\bet^{q-2}\Big(\frac{1}{\|{\bf y}\|}\min\limits_{1\le k\le q-1}{\bf y}'{\bf z}_k^q \Big),\label{maxLECD}\\
A^{q-2}(\Pi({\bf y}))&=\cos^{-1}(t(\Pi({\bf y})))=\ts\cos^{-1}\Big(\frac{1}{\|{\bf y}\|}\min\limits_{1\le k\le q-1}{\bf y}'{\bf z}_k^q \Big).\label{maxLAD}
\end{align} 
Thus,  because $\bet^{q-2}(\cdot)$ and $\cos^{-1}(\cdot)$ are decreasing and $\frac{{\bf y}}{\|{\bf y}\|}$ is a unit vector, we seek a unit vector ${\bf z}\equiv(z_1,\dots,z_q)'\in\~{\cal M}_\le^{q-1}$ that attains the maximum
\begin{equation}\label{Lambdahat}
\^\Lam_q:=\max_{{\bf z}\in\~{\cal M}_\le^{q-1},\,\|{\bf z}\|=1} \; \min_{1\le k\le q-1}{\bf z}'{\bf z}_k^q.
\end{equation}

For $1\le k\le q$ define 
\begin{align}
b_k^q&=\ts\sqrt{\frac{3k(q-k)}{q(q+1)}},\qquad(b_0^q=0),\label{ckdef}\\
\^a_k^q&=b_{k-1}^q-b_k^q,\label{ckdef1}\\
\^{\bf a}^q&=(\^a_1^q,\dots,\^a_q^q)'\label{ckdef2},\\
\^{\bf z}^q&\equiv(\^z_1^q,\dots,\^z_q^q)'=\ts \frac{\^{\bf a}^q}{\|\^{\bf a}^q\|}.\label{zqdef}
\end{align}
Then $\^a_1^q+\cdots+\^a_q^q=0$ so $\^z_1^q+\cdots+\^z_q^q=0$,  and it is straightforward to show that $\^a_1^q<\cdots<\^a_q^q$, so $\^z_1^q<\cdots<\^z_q^q$, hence $\^{\bf z}^q\in\~{\cal M}_\le^{q-1}$. Trivially, $\|\^{\bf z}^q\|=1$.
\vskip4pt

\nid{\bf Proposition 4.1.} The unit vector $\^{\bf z}^q$ uniquely attains  the maximum $\^\Lam_q$. Thus in the original scale,
\begin{equation}\label{yhatq}
\^{\bf y}^q:=\|\bar{\bf y}^q\|\;\^{\bf z}^q=\ts\sqrt{\frac{q(q^2-1)}{12}}\,\frac{\^{\bf a}^q}{\|\^{\bf a}^q\|}
\end{equation}
uniquely minimizes the 
LECD and the LECAD of $\Pi({\bf y})$ for ${\bf y}\in\~{\cal S}_{\|\bar{\bf y}^q\|}^{q-2}\cap\mathbb{R}_\le^q$,
and $\^{\bf y}^q\ne\bar{\bf y}^q$ when $q\ge4$. The minimum LECD and LECAD are
 \begin{align}
 L^{q-2}(\Pi(\^{\bf y}^q))&= \bet^{q-2}\Big(\ts\frac{1}{\|\^{\bf a}^q\|}\sqrt{\frac{3}{q+1}}\,\Big),\label{cosinaoverL}\\
A^{q-2}(\Pi(\^{\bf y}^q))&= \cos^{-1}\Big(\ts\frac{1}{\|\^{\bf a}^q\|}\sqrt{\frac{3}{q+1}}\,\Big).\label{cosinaovera}
 \end{align}

\nid{\bf Proof.} For any unit vector ${\bf z}\equiv(z_1,\dots,z_q)'\in\~{\cal M}_\le^{q-1}$, $z_1+\cdots+z_q=0$, so after some algebra we find that
\begin{align}
{\bf z}'{\bf z}_k^q&=
\ts\sqrt{\frac{q}{k(q-k)}}\;(z_{k+1}+\cdots+z_q),\label{min4}
\end{align}
hence
\begin{equation}\label{min5}
\^\Lam_q=\max_{{\bf z}\in\~{\cal M}_\le^{q-1},\,\|{\bf z}\|=1} \; \min_{1\le k\le q-1}\,\ts\sqrt{\frac{q}{k(q-k)}}\;(z_{k+1}+\cdots+z_q).
\end{equation}
We now show that the maximum in \eqref{min5} is uniquely attained when ${\bf z}=\^{\bf z}^q$.

Because $\^z_{k+1}^q+\cdots+\^z_q^q=\frac{b_k^q}{\|a^q\|}$, 
\begin{equation}\label{min6}
\ts\sqrt{\frac{q}{k(q-k)}}\;(\^z_{k+1}^q+\cdots+\^z_q^q)=\frac{1}{\|\^{\bf a}^q\|}\sqrt{\frac{3}{q+1}}
\end{equation}
for each $k=1,\dots,q-1$. Thus we must show that
\begin{equation}\label{min7}
\min_{1\le k\le q-1}\,\ts\sqrt{\frac{q}{k(q-k)}}\;(z_{k+1}+\cdots+z_q)<\frac{1}{\|\^{\bf a}^q\|}\sqrt{\frac{3}{q+1}}
\end{equation}
for every ${\bf z}\ne \^{\bf z}^q$ such that $z_1+\cdots+z_q=0,\,\|{\bf z}\|=1,\,z_1\le\cdots\le z_q$. Suppose that there is such a ${\bf z}$ that satisfies
\begin{equation}\label{min71}
\min_{1\le k\le q-1}\,\ts\sqrt{\frac{q}{k(q-k)}}\;(z_{k+1}+\cdots+z_q)\ge\frac{1}{\|\^{\bf a}^q\|}\sqrt{\frac{3}{q+1}}.
\end{equation}
Therefore if $1\le k\le q-1$ then
\begin{equation*}
z_{k+1}+\cdots+z_q\ge \ts\frac{b_k^q}{\|\^{\bf a}^q\|}=\^z_{k+1}^q+\cdots+\^z_q^q,
\end{equation*}
with equality for $k=0$, so ${\bf z}$ {\it majorizes} $\^{\bf z}^q$ (Marshall and Olkin [MO]). Because $\|{\bf z}\|^2$ is symmetric and strictly convex in $(z_1,\dots,z_q)$ and ${\bf z}\ne\^{\bf z}^q$, this implies that
\begin{equation}\label{min9}
1=\|{\bf z}\|^2>\|\^{\bf z}^q\|^2=1,
\end{equation}
 a contradiction. Thus the maximum value $\^\Lam^q$ is uniquely achieved when ${\bf z}=\^{\bf z}^q$ as asserted. It is easy to verify that $\^a_1^q,\dots,\^a_q^q$ are not evenly spaced when $q\ge4$, hence $\^{\bf y}^q\ne\bar{\bf y}^q$. Lastly, \eqref{cosinaoverL} and \eqref{cosinaovera} follow from \eqref{min6}.\hfill$\square$

\vskip4pt

The vectors $\^{\bf y}^q$ and $\^{\bf x}^q\equiv\^{\bf y}^q+\frac{q+1}{2}{\bf e}^q$ are called the {\it maximal configurations} in $\~{\cal M}_\le^{q-1}$ and $\mathbb{R}_\le^q$ respectively. It is now obvious to ask whether or not the sequences $\{\^{\bf y}^q\}$ and $\{\^{\bf x}^q\}$ are APF, and if so, are APU. These questions will be answered in Propositions 4.5 and 4.7.

Because the LECD of $\Pi(\^{\bf y}^q)$ given by \eqref{cosinaoverL} depends on $\|\^{\bf a}^q\|$, bounds for $\|\^{\bf a}^q\|$ are needed. Since $\^{\bf y}^q\ne \bar{\bf y}^q$,  necessarily $\|\^{\bf a}^q\|<1$ by the uniqueness of $\^{\bf y}^q$, but sharper bounds will be required.
\vskip4pt

\nid{\bf Lemma 4.2.}
\begin{align}\label{aqbounds}
\ts\sqrt{\frac{3[\log(2q+1)-2]}{2(q+1)}}<\|\^{\bf a}^q\|<\sqrt{\ts \frac{3[2\log(2q-1)+1]}{2(q+1)}}.
\end{align}
\nid Therefore
\begin{equation}\label{aqasympt}
\ts\|\^{\bf a}^q\|=O\Big(\sqrt{\frac{\log q}{q+1}}\,\Big)\ \ {\rm as}\ q\to\infty.
\end{equation}
\nid{\bf Proof.}  For $k=1,\dots,q$ set
\begin{align}
c_k^q&=\ts\big[\sqrt{(k-1(q-k+1)}-\sqrt{k(q-k)}\big]^2,\label{ckdk}\\
d_k^q&=\ts\sqrt{k(k-1)(q-k)(q-k+1)},\nonumber\\
\bar q&=\ts\frac{q+1}{2},\nonumber
\end{align}
then verify that
\begin{equation}\label{ckdk1}
c_k^q=\ts2\big[\frac{q^2-1}{4}-(k-\bar q)^2-d_k^q\big].
\end{equation}
From \eqref{ckdef}-\eqref{ckdef1} and \eqref{ckdk}-\eqref{ckdk1} we find that
\begin{align}
\|\^{\bf a}^q\|^2
&=\ts\frac{3}{q(q+1)}\sum\nolimits_{k=1}^qc_k^q\label{sumsqs}\\
 &=\ts(q-1)-\ts\frac{6}{q(q+1)}\sum\limits_{k=1}^qd_k^q.\nonumber
 \end{align}
 
For the upper bound, use the harmonic mean-geometric mean inequality:
\begin{align}
\|\^{\bf a}^q\|^2
<&\;(q-1)-\ts\frac{6}{q(q+1)}\sum\nolimits_{k=1}^q
 \frac{k(k-1)(q-k)(q-k+1)}{(k-\frac{1}{2})(q-k+\frac{1}{2})}\nonumber\\
 =&\;(q-1)-\ts\frac{6}{q(q+1)}\sum\nolimits_{k=1}^q
 \frac{[(k-\frac{1}{2})^2-\frac{1}{4}][(q-k+\frac{1}{2})^2-\frac{1}{4}]}{(k-\frac{1}{2})(q-k+\frac{1}{2})}\nonumber\\
 <&\;(q-1)-\;\ts\frac{6}{q(q+1)}\sum\nolimits_{k=1}^q\Big\{
 (k-\frac{1}{2})(q-k+\frac{1}{2})-\frac{k-\frac{1}{2}}{4(q-k+\frac{1}{2})}-\frac{q-k+\frac{1}{2}}{4(k-\frac{1}{2})}
 \Big\}\nonumber\\
 =&\;(q-1)-\ts\frac{6}{q(q+1)}\Big\{\sum\nolimits_{k=1}^q
 (k-\frac{1}{2})(q-k+\frac{1}{2})-\sum\nolimits_{k=1}^q\frac{q-k+\frac{1}{2}}{2(k-\frac{1}{2})}\Big\}\nonumber\\
 =&\;(q-1)-\ts\frac{6}{q(q+1)}\Big\{\sum\nolimits_{k=1}^q
 (k-\frac{1}{2})(q-k+\frac{1}{2})-\frac{q}{2}\sum\nolimits_{k=1}^q\frac{1}{k-\frac{1}{2}}+\frac{q}{2}\Big\}\nonumber\\
 =&\;\ts\frac{3}{q+1}\Big\{\sum\nolimits_{k=1}^q\frac{1}{k-\frac{1}{2}}-\frac{3}{2}\Big\}\nonumber\\
 <&\,\ts \frac{3[2\log(2q-1)+1]}{2(q+1)};\label{upperbd10}
\end{align}
the final inequality follows from (7) of Qi and Guo [QG].

Similarly, the geometric mean-arithmetic mean inequality yields the non-logarithmic lower bound $\frac{3(3q-2)}{2q(q+1)}$. However, the asserted logarithmic lower bound, which is sharper, can be obtained as follows. We will show that 
\begin{align}
c_k^q\equiv\big[\sqrt{k(q-k)}-\sqrt{(k-1)(q-k+1)}\,\big]^2\ge\ts\frac{(k-\bar q)^2}{(k-\frac{1}{2})(q-k+\frac{1}{2})},\label{sumsqs3}
\end{align}
for $k=1,\dots,q$, where $\bar q=\frac{q+1}{2}$. Thus from \eqref{sumsqs}, 
\begin{align}
\|\^{\bf a}^q\|^2
 \ge&\;\ts\frac{3}{q(q+1)}\sum_{k=1}^q\frac{(k-\bar q)^2}{(k-\frac{1}{2})(q-k+\frac{1}{2})}\nonumber\\
 =&\;\ts\frac{3}{q^2(q+1)}\sum_{k=1}^q\Big[\frac{(k-\bar q)^2}{k-\frac{1}{2}}+\frac{(k-\bar q)^2}{q-k+\frac{1}{2}}\Big]\nonumber\\
  =&\;\ts\frac{3}{q^2(q+1)}\sum_{k=1}^q\Big[\frac{(k-\frac{1}{2})^2-2(k-\frac{1}{2})(\frac{q}{2})+(\frac{q}{2})^2}{k-\frac{1}{2}}+\frac{(q-k+\frac{1}{2})^2-2(q-k+\frac{1}{2})(\frac{q}{2})+(\frac{q}{2})^2}{q-k+\frac{1}{2}}\Big]\nonumber\\
 =&\;\ts\frac{3}{q+1}\Big[\sum_{k=1}^q\frac{1}{4}\big(\frac{1}{k-\frac{1}{2}}+\frac{1}{q-k+\frac{1}{2}}\big)-1\Big]\nonumber\\
  =&\;\ts\frac{3}{q+1}\big[\sum_{k=1}^q\frac{1}{2k-1}-1\big]\nonumber\\
  >&\;\ts\frac{3[\log(2q+1)-2]}{2(q+1)},\label{sumsqs4}
 \end{align} 
 where the inequality used in \eqref{sumsqs4} also follows from (7) of [QG].
 
To establish \eqref{sumsqs3}, rewrite it in the equivalent form
\begin{align}
\big(\sqrt{(\bar q+u)(\bar q-u-1)}-\sqrt{(\bar q+u-1)(\bar q-u)}\,\big)^2\ge\ts\frac{u^2}{ \~q^2-u^2},\label{sumsqs7}
\end{align}
where $u\equiv k-\bar q\in\{-\frac{q-1}{2},\dots,\frac{q-1}{2}\}$ and $\~q=\bar q-\t12=\frac{q}{2}$. Now set $v=\frac{u}{\~q}$, so $|v|\le\frac{q-1}{q}<1$. Then \eqref{sumsqs7} can be written in the equivalent forms
\begin{align}
\big(\sqrt{(\bar q+\~q v)(\bar q-\~q v-1)}-\sqrt{(\bar q+\~q v-1)(\bar q-\~q v)}\,\big)^2&\ge\ts\frac{v^2}{ 1-v^2},\nonumber\\
\big(\sqrt{(\bar q^2-\~q^2 v^2)-(\bar q+\~q v)}-\sqrt{(\bar q^2-\~q^2 v^2)-(\bar q-\~q v)}\,\big)^2&\ge\ts\frac{v^2}{ 1-v^2},\nonumber\\
2\mu(v)-2\sqrt{(\mu(v)-\~q v)(\mu(v)+\~q v)}&\ge\ts\frac{v^2}{ 1-v^2},\nonumber\\
2\mu(v)-2\sqrt{\mu(v)^2-\~q^2 v^2}&\ge\ts\frac{v^2}{ 1-v^2}\label{sumsqs8}
\end{align}
where
\begin{equation*}
\mu(v)=\bar q^2-\~q^2v^2-\bar q=\ts\frac{1}{4}[q^2(1-v^2)-1].
\end{equation*}
It will be shown that for $|v|\le\frac{q-1}{q}$, 
\begin{equation}\label{oktosquare}
2\mu(v)-\ts\frac{v^2}{ 1-v^2}\ge0,
\end{equation}
so \eqref{sumsqs8} is equivalent to each of the following inequalities:
\begin{align}
[2\mu(v)-\ts\frac{v^2}{ 1-v^2}]^2&\ge4\mu(v)^2-q^2v^2,\nonumber\\
q^2v^2-4\mu(v)\ts\frac{v^2}{ 1-v^2}+\ts\frac{v^4}{ (1-v^2)^2}&\ge0,\nonumber\\
q^2v^2-[q^2(1-v^2)-1]\ts\frac{v^2}{ 1-v^2}+\ts\frac{v^4}{ (1-v^2)^2}&\ge0,\nonumber\\
q^2v^2(1-v^2)-[q^2(1-v^2)-1]v^2+\ts\frac{v^4}{ 1-v^2}&\ge0,\nonumber\\
v^2(1+\ts\frac{v^2}{ 1-v^2})&\ge0,\nonumber
\end{align}
which clearly is true. Thus \eqref{sumsqs3} will be established once \eqref{oktosquare} is verified.  

For this set $x=v^2$, so \eqref{oktosquare} can be expressed equivalently as
\begin{align*}
h(x)\equiv(1-x)[q^2(1-x)-1]-2x&\ge0,
\end{align*}
where $0\le x\le(\frac{q-1}{q})^2$. The quadratic function $h(x)$ satisfies
\begin{equation*}
\ts h(0)=q^2-1>h\big[\big(\frac{q-1}{q}\big)^2\big]=2-\frac{2}{q}>0>h(1)=-2,
\end{equation*}
hence $h(x)>0$ for $0\le x\le(\frac{q-1}{q})^2$, as required.
\hfill$\square$
\newpage

\nid{\bf Proposition 4.5.} The maximal configurations $\{\^{\bf y}^q\}$ and $\{\^{\bf x}^q\}$ are APF.
\vskip2pt

\nid{\bf Proof.} Set $t^q=\ts\frac{1}{\|\^{\bf a}^q\|}\sqrt{\frac{3}{q+1}}$, so that \eqref{aqbounds} yields
\begin{equation*}
\ts\sqrt{ \frac{q}{2\log(2q-1)+1} }<t^q\sqrt{q}<\sqrt{ \frac{2q}{\log(2q+1)-2} }.
\end{equation*}
Then by Lemma 2.1 with $\lam=\infty$, 
\begin{equation}\label{Lhatlimit}
\lim_{q\to\infty}L^{q-2}(\Pi(\^{\bf y}^q))=0,
\end{equation}
hence $\{\^{\bf y}^q\}$ (and $\{\^{\bf x}^q\}$) is APF.\hfill$\square$

\vskip4pt

It follows from \eqref{coordNormal}, \eqref{normyq}, and \eqref{Uniformlim} that for each fixed $i\ge1$,
\begin{align}
 W_q:=\ts\frac{12}{q^2-1}\big[(\~{\bf U}_{\|\bar{\bf y}^q\|}^{q-2})_i\big]^2&\tod \chi_1^2,\label{chisquarelim}\\
 \bar W_q:=\ts\frac{12}{q^2-1}\big[(\~{\bf U}_{\bar{\bf y}^q}^{q-2})_i\big]^2&\tod3\,\mathrm{Beta}(\t12,1),\label{Betalim}
\end{align}
as $q\to\infty$. The bounds for $\|\^{\bf a}^q\|$ in \eqref{aqbounds} yield a corresponding result for the maximal configuration:
\vskip4pt

\nid{\bf Proposition 4.6.} For each fixed $i\ge1$,
\begin{align}
\ts\frac{Z_q}{\log(2q-1)+2} <_{\mathrm{st}}\^ W_q:=\frac{12}{q^2-1}\big[(\~{\bf U}_{\^{\bf y}^q}^{q-2})_i\big]^2<_{\mathrm{st}}\frac{2Z_q+1}{\log(2q+1)-2},\label{Maxconfiglim1}
\end{align}
where $\{Z_q\}$ is a sequence of positive random variables such that 
\begin{align}
Z_q&\ts\tod F_{1,2}\label{Yqlim}
\end{align}
as $q\to\infty$. Here $<_{\mathrm{st}}$ denotes stochastic ordering and $F_{1,2}$ denotes the F distribution with 1 and 2 degrees of freedom. Therefore
\begin{align}
\ts\^W_q&\ts=O_p\Big(\frac{F_{1,2}}{\log q}\Big)\top0,\label{hatWlim}\\
\ts\sqrt{\frac{12}{q^2-1}}(\~{\bf U}_{\^{\bf y}^q}^{q-2})_i &\ts=O_p\big(\frac{t_2}{\sqrt{\log q}}\big)\top0,\label{hatWlim2}
\end{align}
where $t_2$ denotes Student's t-distribution with 2 degrees of freedom. 
\vskip2pt
\nid{\bf Proof.} From \eqref{ckdef}-\eqref{yhatq} and \eqref{ckdk}, $\^W_q$ is uniformly distributed over the set
\begin{align}
\ts\big\{\frac{12}{q^2-1}\frac{q(q^2-1)}{12}\frac{3}{\|\^{\bf a}^q\|^2q(q+1)}\,c_k^q\Bigm| k=1,\dots,q\big\}
=\big\{\frac{3c_k^q}{(q+1)\|\^{\bf a}^q\|^2}\bigm| k=1,\dots,q\big\},\label{hatWunifdist}
\end{align}
so by \eqref{sumsqs4}, $\^W_q$ is stochastically smaller than the uniform distribution on 
\begin{align*}
\ts\Big\{\frac{2c_k^q}{\log(2q+1)-2}\Bigm| k=1,\dots,q\Big\}=\Big\{\frac{4[\frac{q^2-1}{4}-(k-\bar q)^2-d_k^q]}{\log(2q+1)-2}\Bigm| k=1,\dots,q\Big\}.
\end{align*}
Now apply the harmonic mean-geometric mean inequality to $d_k^q$ to obtain
\begin{align*}
\ts\frac{q^2-1}{4}-(k-\bar q)^2-d_k^q&<\ts\frac{q^2-1}{4}-(k-\bar q)^2- \frac{k(k-1)(q-k)(q-k+1)}{(k-\frac{1}{2})(q-k+\frac{1}{2})}\\
 &=\ts\frac{q^2-1}{4}-(k-\bar q)^2-\frac{[(k-\frac{1}{2})^2-\frac{1}{4}][(q-k+\frac{1}{2})^2-\frac{1}{4}]}{(k-\frac{1}{2})(q-k+\frac{1}{2})}\\
 &=\ts\frac{q^2-1}{4}-(k-\bar q)^2-(k-\frac{1}{2})(q-k+\frac{1}{2})\\
  &\ts\ \ \ +\frac{q-k+\frac{1}{2}}{4(k-\frac{1}{2})} +\frac{k-\frac{1}{2}}{4(q-k+\frac{1}{2})}
  -\frac{1}{16}\frac{1}{(k-\frac{1}{2})(q-k+\frac{1}{2})}\\
&=\ts\frac{1}{4}\Big[\frac{q^2-\frac{1}{4}}{(k-\frac{1}{2})(q-k+\frac{1}{2})}-3\Big]\\
&=\ts\frac{1}{4}\Big[\frac{q^2-\frac{1}{4}} {\frac{q^2}{4}-(k-\bar q)^2}-3\Big]\\
&<\ts\frac{1} {1-4\big(\frac{k-\bar q}{q}\big)^2}-\frac{3}{4},
\end{align*}
where we have twice used the relation
\begin{equation}\label{relation}
\ts(k-\bar q)^2+(k-\frac{1}{2})(q-k+\frac{1}{2})=\frac{q^2}{4}.
\end{equation}
Therefore $\^W_q$ is stochastically smaller than 
\begin{equation*}
\ts\frac{4}{\log(2q+1)-2}\big(\frac{1} {1-4V_q^2}-\frac{3}{4}\big)\equiv\frac{4Y_q+1}{\log(2q+1)-2},
\end{equation*}
where 
\begin{align*}
V_q&\ts\eqd\mathrm{Uniform}\big\{\frac{k-\bar q}{q}\bigm| k=1,\dots,q\big\},\\
Y_q&\ts=\frac{4V_q^2} {1-4V_q^2}.
\end{align*}
Because $\frac{k-\bar q}{q}=\frac{2k-q-1}{2q}$, clearly
\begin{align*}
V_q&\tod V\eqd\rm{Uniform}(-\t12,\t12),\\
4V_q^2&\tod 4V^2\eqd\mathrm{Beta}(\t12,1),
\end{align*}
as $q\to\infty$, from which it follows that $Y_q\tod \t12F_{1,2}$. Now set $Z_q=2Y_q$.

Similarly from \eqref{upperbd10}, \eqref{sumsqs3}, \eqref{hatWunifdist}, and \eqref{relation}, $\^W_q$ is stochastically larger than the uniform distribution on 
\begin{equation*}
\ts\Big\{\frac{2}{\log(2q+1)-2}\Big[\frac{(k-\bar q)^2}{\frac{q^2}{4}-(k-\bar q)^2}\Big]\Bigm| k=1,\dots,q\Big\},
\end{equation*}
so $\^W_q$ is stochastically larger than
\begin{equation*}
\ts\frac{2}{\log(2q+1)-2}\big(\frac{4V_q^2} {1-4V_q^2}\big)\equiv\frac{Z_q}{\log(2q+1)-2},
\end{equation*}
as asserted. \hfill$\square$
\vskip4pt


\nid{\bf Proposition 4.7.} The sequences of maximal configurations $\{\^{\bf y}^q\}$ and $\{\^{\bf x}^q\}$ are not APU.
\vskip2pt

\nid{\bf Proof.} It follows from \eqref{hatWlim2}
that for any fixed $i$,
\begin{equation}\label{Uilimit}
\ts\frac{\sqrt{q}}{\|\^{\bf y}^q\|}(\~{\bf U}_{\^{\bf y}^q}^{q-2})_i= \sqrt{\frac{12}{q^2-1}}(\~{\bf U}_{\^{\bf y}^q}^{q-2})_i=O_p\big(\frac{1}{\sqrt{\log q}}\big)\top0,
\end{equation}
hence by Proposition 2.10 $\{\^{\bf y}^q\}$ and $\{\^{\bf x}^q\}$ cannot be APU.\hfill$\square$
\vskip4pt

\vskip6pt

\nid{\bf 5. The normal configuration.}

\nid  The sequence $\{\^{\bf y}^q\}$, like $\{\bar{\bf y}^q\}$, fails to satisfy the necessary condition \eqref{necCondAPU} for APU, yet $\{\^{\bf y}^q\}$ uniquely minimizes the LECD and LECAD, so it seems reasonable to conjecture that no APU sequence exists. However, it is easy to find a sequence $\{\u{\bf y}^q\in\~{\cal M}_\le^{q-1}\}$ that does satisfy \eqref{necCondAPU}.
Define
\begin{equation}\label{acupq}
\ts \u{\bf a}^q\equiv(\u a_1^q,\dots,\u a_q^q)=\big(\Phi^{-1}(\frac{1}{q+1}),\Phi^{-1}(\frac{2}{q+1}),\dots,\Phi^{-1}(\frac{q}{q+1})\big)',
\end{equation}
the $\frac{k}{q+1}$-quantiles of the $N(0,1)$ distribution, then in the original scale let
\begin{equation}\label{yuq}
\ts\u{\bf y}^q=\|\bar{\bf y}^q\|\frac{\u{\bf a}^q}{\|\u{\bf a}^q\|}.
\end{equation}
Clearly $\u a_1^q<\cdots<\u a_q^q$ while $\u a_1^q+\cdots+\u a_q^q=0$ by the symmetry of $N(0,1)$, hence $\u{\bf y}^q\in\~{\cal M}_\le^{q-1}$. The vector $\u{\bf y}^q$ is called the {\it normal configuration}. 

For each $i=1,\dots,q$, $(\~{\bf U}_{\u{\bf a}^q}^{q-2})_i\eqd\Phi^{-1}(U_q)$, where 
\begin{equation*}
\ts U_q\eqd\mathrm{Uniform}\big(\big\{\Phi^{-1}(\frac{1}{q+1}),\dots,\Phi^{-1}(\frac{q}{q+1})\big\}\big)\tod\mathrm{Uniform}(0,1),
\end{equation*}
hence $(\~{\bf U}_{\u{\bf a}^q}^{q-2})_i\tod N(0,1)$ as $q\to\infty$. Furthermore,
\begin{equation*}
\ts\frac{\|\u{\bf a}^q\|^2}{q+1}+\frac{1}{q+1}[\Phi^{-1}(\frac{q}{q+1})\big]^2
\equiv\frac{1}{q+1}\sum_{k=1}^q\big[\Phi^{-1}(\frac{k}{q+1})\big]^2+\frac{1}{q+1}[\Phi^{-1}(\frac{q}{q+1})\big]^2
\end{equation*}
is an approximating Riemann sum for
\begin{equation*}
\ts\int_0^1[\Phi^{-1}(u)]^2du=\int_{-\infty}^\infty x^2\phi(x)dx=1,
\end{equation*}
while
\begin{equation}\label{Phiinvapprox}
\ts\Phi^{-1}(\frac{q}{q+1})=\sqrt{2\log(q+1)}(1+o(1))
\end{equation}
as $q\to\infty$ (e.g. Fung and Seneta [FS] p.1092), hence
\begin{align}
\ts\frac{\|\u{\bf a}^q\|^2}{q+1}&=\ts 1-\frac{2\log(q+1)}{q+1}+o(1),\label{anorm1}\\
\|\u{\bf a}^q\|&\sim\sqrt{q+1}.\label{anorm}
\end{align}
Therefore $\{\u{\bf y}^q\}$ satisfies \eqref{necCondAPU}:
\begin{equation}\label{necCondAPU1}
\ts\frac{\sqrt{q}}{\|\u{\bf y}^q\|}(\~{\bf U}_{\u{\bf y}^q}^{q-2})_i=\frac{\sqrt{q}}{\|\u{\bf a}^q\|}(\~{\bf U}_{\u{\bf a}^q}^{q-2})_i\tod N(0,1)
\end{equation}
as $q\to\infty$. However, it is now shown that the LECD of $\{\u{\bf y}^q\}$, necessarily greater than that of $\{\^{\bf y}^q\}$, does not approach 0.
\vskip4pt

\nid{\bf Proposition 5.1.} $\{\u{\bf y}^q\}$ is not APF, hence is not APU.
\vskip2pt

\nid{\bf Proof.} By \eqref{LECDdef.2}-\eqref{LleD1} and \eqref{acupq}-\eqref{yuq}, 
\begin{align}
L^{q-2}\big(\Pi(\u{\bf y}^q)\big)&=\ts\bet^{q-2}(\u t^q),\label{Ltuq}\\
\ts\u t^q:&=\ts\frac{1}{\|\u{\bf y}^q\|}\min\limits_{1\le k\le q-1} (\u{\bf y}^q)'{\bf z}_k^q\label{tuq}\\
&\le\ts\ts\frac{1}{\|\u{\bf y}^q\|}(\u{\bf y}^q)'{\bf z}_{q-1}^q\nonumber\\
&=\ts\frac{1}{\|\u{\bf a}^q\|}(\u{\bf a}^q)'{\bf z}_{q-1}^q\label{tuqstar}\\
&=\ts\frac{1}{\|\u{\bf a}^q\|}\sqrt{\frac{q}{q-1}}\Phi^{-1}\big(\frac{q}{q+1}\big)\label{tuqstar1}\\
&<\ts\frac{1}{\|\u{\bf a}^q\|}\sqrt{\frac{q}{q-1}}\frac{\phi(\Phi^{-1}(\frac{q}{q+1}))}{1-\Phi(\Phi^{-1}(\frac{q}{q+1}))}\nonumber\\
&=\ts\frac{q+1}{\sqrt{2\pi}\|\u{\bf a}^q\|}\sqrt{\frac{q}{q-1}}e^{-\frac{1}{2}[\Phi^{-1}(\frac{q}{q+1})]^2}.\nonumber
\end{align}
It follows from Fung and Seneta [FS] p.1092 that
\begin{equation}\label{FungSen}
\ts\Phi^{-1}\big(\frac{q}{q+1}\big)=\sqrt{2\log\big((q+1)\sqrt{4\pi\log(q+1)}\,\big)}\bigg(1+\Del_q\bigg)
\end{equation}
where $\Del_q=O\Big(\frac{\log(\log(q+1))}{(\log(q+1))^2}\Big)$, hence
\begin{align*}
\ts e^{-\frac{1}{2}[\Phi^{-1}(\frac{q}{q+1})]^2}
&=\ts\frac{1}{q+1}\frac{1}{\sqrt{4\pi\log(q+1)}}\big(\frac{1}{q+1}\big)^{\Del_q+\Del_q^2}\big(\frac{1}{\sqrt{4\pi\log(q+1)}}\big)^{\Del_q+\Del_q^2}\\
&=\ts\frac{1}{q+1} o(1) (1+o(1)) (1+o(1))\\
&=\ts \frac{1}{q+1} o(1).
\end{align*}
Therefore by \eqref{anorm},
\begin{align}
\u t^q\sqrt{q}<\ts\frac{1}{\sqrt{2\pi}\|\u{\bf a}^q\|}\frac{q}{\sqrt{q-1}}o(1)=o(1),\label{tqo1}
\end{align}
hence $\lim_{q\to\infty}\u t^q\sqrt{q}=0$, so
\begin{equation}\label{limLycup}
\lim_{q\to\infty}L^{q-2}\big(\Pi(\u{\bf y}^q)\big)=\t12
\end{equation}
by Lemma 2.1 with $\lam=0$. This completes the proof.\hfill$\square$
\vskip4pt

\nid{\bf Remark 5.2.} It should be noted that the convergences in \eqref{Uilimit} and \eqref{limLycup} occur at very slow, sub-logarithmic rates.\hfill$\square$
\vskip4pt

\nid{\bf Proposition 5.3.} $\{\u{\bf y}^q\}$ is APE.
\vskip2pt

\nid{\bf Proof.} The proof is similar to that of Proposition 3.2. Again define $\bar{\bf z}_k^q=\|\bar{\bf y}^q\|{\bf z}_k^q$, where ${\bf z}_k^q$ is the unit vector in \eqref{ak}, and define
\begin{equation*}
\ts\u s^q=\frac{1}{\|\u{\bf a}^q\|}(\u{\bf a}^q)'{\bf z}_{q-1}^q,
\end{equation*}
cf. \eqref{tuqstar}. As in \eqref{tuqstar}-\eqref{tqo1},
\begin{align}
\u s^q\sqrt{q}<\ts\frac{1}{\sqrt{2\pi}\|\u{\bf a}^q\|}\frac{q}{\sqrt{q-1}}o(1)=o(1),
\end{align}
hence from \eqref{Betabd}-\eqref{betdef} and Lemma 2.1 with $\lam=0$,
\begin{align}
\lim_{q\to\infty}\~{\bf U}_{\|\bar{\bf y}^q\|}^{q-2}\big(C\big(\bar{\bf z}_{q-1}^q;\u s^q\big)\big)=\lim_{q\to\infty}\bet^{q-2}(\u s^q)=\t12.\label{Luyq}
\end{align}
Furthermore, from \eqref{empcapcond} and the Rearrangement Inequality,
\begin{align*}
\Pi(\u{\bf y}^q)\cap C\big(\bar{\bf z}_{q-1}^q;\u s^q\big)=\nul&\iff\max_{P\in{\cal P}^q}(P\u{\bf y}^q)'\bar{\bf z}_{q-1}^q\le\|\bar{\bf y}^q\|^2\u s^q\\
&\iff(\u{\bf y}^q)'\bar{\bf z}_{q-1}^q\le\|\bar{\bf y}^q\|^2\u s^q\\
&\iff\ts\frac{1}{\|\u{\bf a}^q\|}(\u{\bf a}^q)'{\bf z}_{q-1}^q\le \u s^q,
\end{align*}
hence
$C\big(\bar{\bf z}_{q-1}^q;\u s^q\big)$ is an empty spherical cap for $\Pi(\u{\bf y}^q)$. 

Because $P\Pi(\u{\bf y}^q)=\Pi(\u{\bf y}^q)$ for all $P\in{\cal P}^q$, each $C\big(P\bar{\bf z}_{q-1}^q;\u s^q\big)$ is an empty spherical cap for $\Pi(\u{\bf y}^q)$ in $\~{\cal S}_{\|\bar{\bf y}^q\|}^{q-2}$; there are $q!$ such congruent caps. However
\begin{align*}
\{P\bar{\bf z}_{q-1}^q\mid P\in{\cal P}^q\}&=\{\bar{\bf f}_1^q,\dots,\bar{\bf f}_q^q\},
\end{align*}
where $\bar{\bf f}_i^q=\|\bar{\bf y}^q\|{\bf f}_i^q$, so these $q!$ empty caps reduce to $q$ congruent ones, namely
\begin{equation*}
\ts\big\{C\big(\bar{\bf f}_i^q;\u s^q\big)\bigm| i=1,\dots q\big\}.
\end{equation*}
By \eqref{Luyq} each of these congruent caps remains nonnegligible as $q\to\infty$, so to show that $\{\bar{\bf y}^q\}$ is APE it suffices to show that
\begin{equation}\label{limuUps}
\lim_{q\to\infty}\~{\bf U}_{\|\bar{\bf y}^q\|}^{q-2}\big(\u\Ups^q\big)=1,
\end{equation}
where
\begin{align*}
\u\Ups^q&=\ts\bigcup\nolimits_{i=1}^q \Big[ C\big(\bar{\bf f}_i^q;\u s^q\big)\Big].
\end{align*}

Clearly
\begin{equation}\label{uUpsinclusion}
\ts\~{\cal S}_{\|\bar{\bf y}^q\|}^{q-2}\bigcap\big((\u\Ups_q)^c\big)\subseteq\~{\cal S}_{\|\bar{\bf y}^q\|}^{q-2}\bigcap\big(\cap_{i=2}^q \u H_i^q\big),
\end{equation}
where $\u H_i^q$ is the halfspace
\begin{align*}
\u H_i^q&:=\ts\big\{{\bf v}\in\~{\cal M}^{q-1}\bigm| {\bf v}'\bar{\bf f}_i^q\le\|\bar{\bf y}^q\|^2\,\u s^q\big\}.
\end{align*}
As in the proof of Proposition 3.2, $({\bf f}_i^q)'{\bf f}_j^q=-\frac{1}{q-1}<0$ if $i\ne j$ so 
\begin{equation}\label{aaa1}
\~{\bf U}_{\|\bar{\bf y}^q\|}^{q-2}\big(\cap_{i=2}^q \u H_i^q\big)\le\~{\bf U}_{\|\bar{\bf y}^q\|}^{q-2}\big(\cap_{i=2}^q\u K_i^q\big),
\end{equation}
where $\u K_i^q$ is the halfspace
\begin{equation*}
\ts\u K_i^q:=\big\{{\bf v}\in\~{\cal M}^{q-1}\bigm| {\bf v}'\gam_i^q\le \|\bar{\bf y}^q\|\,\u s^q \big\}.
\end{equation*}
Again apply Proposition A.1 in the Appendix and the orthogonal invariance of $\~{\bf U}_{\|\bar{\bf y}^q\|}^{q-2}$ to obtain 
\begin{align}
\~{\bf U}_{\|\bar{\bf y}^q\|}^{q-2}\big(\cap_{i=2}^q\u K_i^q\big)&\le\ts\prod_{i=2}^q\~{\bf U}_{\|\bar{\bf y}^q\|}^{q-2}(\u K_i^q)\label{uKiqineq1}\\
 &=\ts \big[\~{\bf U}_{\|\bar{\bf y}^q\|}^{q-2}(\u K_i^q)\big]^{q-1}\nonumber\\
  &=\ts\big[1-\bet^{q-2}\big(\u s^q\big)\big]^{q-1}.\label{uKiqineq2}
\end{align}
Therefore by \eqref{Luyq},
\begin{equation}\label{limsupuK}
\ts\limsup\limits_{q\to\infty} \big[\~{\bf U}_{\|\bar{\bf y}^q\|}^{q-2}\big(\cap_{i=2}^q\u K_i^q\big)\big]^{\frac{1}{q-1}}
\le\t12, 
\end{equation}
hence by \eqref{uUpsinclusion}-\eqref{limsupuK},
\begin{equation}\label{Ugeometric1}
\~{\bf U}_{\|\bar{\bf y}^q\|}^{q-2}\big((\Ups^q)^c\big)\le\~{\bf U}_{\|\bar{\bf y}^q\|}^{q-2}\big(\cap_{i=2}^q\u K_i^q\big)\le(.51)^{q-1}
\end{equation}
for sufficiently large $q$. Thus \eqref{limuUps} holds, in fact $\~{\bf U}_{\|\bar{\bf y}^q\|}^{q-2}((\u\Ups^q)^c)\to0$ at a geometric rate, hence $\{\u{\bf y}^q\}$ is APE as asserted.\hfill$\square$
\vskip6pt

\nid{\bf 6. Comparisons among the distributions.} 

\nid Based on the results in \S3-5, comparisons among the three uniform distributions $\~{\bf U}_{\bar{\bf y}^q}^{q-2}$, $\~{\bf U}_{\^{\bf y}^q}^{q-2}$, $\~{\bf U}_{\u{\bf y}^q}^{q-2}$ on permutations and the uniform distribution $\~{\bf U}_{\|\bar{\bf y}^q\|}^{q-2}$ on the sphere $\~{\cal S}_{\|\bar{\bf y}^q\|}^{q-2}$ are now summarized.

The LECDs of $\Pi(\bar{\bf y}^q)$, $\Pi(\^{\bf y}^q)$, and $\Pi(\u{\bf y}^q)$ are as follows:
\begin{align}
L^{q-2}(\Pi(\bar{\bf y}^q))&= \bet^{q-2}\Big(\ts\sqrt{\frac{3}{q+1}}\,\Big);\label{E}\\
L^{q-2}(\Pi(\^{\bf y}^q))&= \bet^{q-2}\Big(\ts\frac{1}{\|\^{\bf a}^q\|}\sqrt{\frac{3}{q+1}}\,\Big);\label{F}\\
L^{q-2}\big(\Pi(\u{\bf y}^q)\big)&=\ts\bet^{q-2}(\u t^q).\label{G}
\end{align}
Here $\|\^{\bf a}^q\|$ is given by \eqref{ckdef2} and approximated in \eqref{aqbounds}, while $\u t^q$ is given by \eqref{tuq} and bounded above by \eqref{tuqstar1} together with \eqref{anorm1}. Some explicit bounds and asymptotic comparisons among these LECDs are collected here.

First, from \eqref{bd108} and \eqref{ineq5},
\begin{align}
\t12-\sqrt{\frac{q-2}{q-4}}\Big[\Phi\Big(\sqrt{\frac{3(q-4)}{q+1}}\Big)-\t12\Big]<L^{q-2}(\Pi(\bar{\bf y}^q))
<\ts\Big[\frac{q-2}{q+1}\Big]^{\frac{q-2}{2}}\sqrt{\frac{q+1}{6\pi(q-2)}}.
\label{bounds23}
\end{align}
Asymptotically,
\begin{equation}\label{Lbarlimit1}
\lim_{q\to\infty}L^{q-2}(\Pi(\bar{\bf y}^q))=1-\Phi(\sqrt{3})\approx0416.
\end{equation}

Second, from \eqref{aqbounds},
\begin{equation*}
\ts\bet^{q-2}\Big(\sqrt{\frac{2}{\log(2q+1)-2}}\Big)<L^{q-2}(\Pi(\^{\bf y}^q))<\bet^{q-2}\Big(\sqrt{\frac{2}{2\log(2q-1)+1}}\Big),
\end{equation*}
which, combined with \eqref{bd108} and \eqref{ineq5}, yields the explicit bounds
\begin{align}
\t12-\sqrt{\frac{q-2}{q-4}}\Big[\Phi\Big(\sqrt{\frac{2(q-4)}{\log(2q+1)-2}}\Big)-\t12\Big]&<L^{q-2}(\Pi(\^{\bf y}^q))\label{bounds21}\\
&<\ts\Big[\frac{2\log(2q-1)-1}{2\log(2q-1)+1}\Big]^{\frac{q-2}{2}}\sqrt{\frac{2\log(2q-1)+1}{4\pi(q-2)}}.\label{bounds22}
\end{align}
Asymptotically,
\begin{equation}\label{Lhatlimit}
\lim_{q\to\infty}L^{q-2}(\Pi(\^{\bf y}^q))=0.
\end{equation}

Third, from \eqref{ineq5} and \eqref{tuqstar1},
\begin{equation*}
\ts\t12-\sqrt{\frac{q-2}{q-4}}\Big[\Phi\Big(\frac{1}{\|\u{\bf a}^q\|}\sqrt{\frac{q(q-4)}{q-1}}\Phi^{-1}\big(\frac{q}{q+1}\big)\Big)-\t12\Big]<L^{q-2}(\Pi(\u{\bf y}^q)).
\end{equation*}
Asymptotically,
\begin{equation}\label{limLycup5}
\lim_{q\to\infty}L^{q-2}\big(\Pi(\u{\bf y}^q)\big)=\t12.
\end{equation}

The LECADs of $\Pi(\bar{\bf y}^q)$, $\Pi(\^{\bf y}^q)$, and $\Pi(\u{\bf y}^q)$ are as follows:
 \begin{align}
A^{q-2}(\Pi(\bar{\bf y}^q))&= \cos^{-1}\Big(\ts\sqrt{\frac{3}{q+1}}\,\Big);\label{H}\\
A^{q-2}(\Pi(\^{\bf y}^q))&= \cos^{-1}\Big(\ts\frac{1}{\|\^{\bf a}^q\|}\sqrt{\frac{3}{q+1}}\,\Big);\label{I}\\
A^{q-2}\big(\Pi(\u{\bf y}^q)\big)&=\ts\cos^{-1}(\u t^q).\label{J}
\end{align}
These yield some explicit expressions and bounds for the LECADs:
\begin{align}
 &A^{q-2}(\Pi(\bar{\bf y}^q))= \cos^{-1}\Big(\ts\sqrt{\frac{3}{q+1}}\,\Big);\label{bounds23A}\\
 \cos^{-1}\Big(\ts\sqrt{\frac{2}{\log(2q+1)-2}}\,\Big)<&A^{q-2}(\Pi(\^{\bf y}^q))<\ts \cos^{-1}\Big(\ts\sqrt{\frac{2}{2\log(2q-1)+1}}\,\Big);\label{bounds22A}\\
  \cos^{-1}\Big(\ts\frac{1}{\|\u{\bf a}^q\|}\sqrt{\frac{q}{q-1}}\Phi^{-1}\big(\frac{q}{q+1}\big)\,\Big)<&A^{q-2}(\Pi(\u{\bf y}^q)).\label{bounds24A}
\end{align}
Asymptotic comparisons among the LECADs are extremely simple:
\vskip4pt

\nid{\bf Proposition 6.1.} For any sequence of nonzero vectors $\{{\bf y}^q\in \~{\cal M}_\le^{q-1}\}$,
\begin{align}
&\lim_{q\to\infty}A^{q-2}(\Pi({\bf y}^q))=\cos^{-1}(0)=\ts\frac{\pi}{2},\label{Ahatlimit4}
\end{align}
that is, the largest empty cap for $\Pi({\bf y}^q)$ approaches a hemisphere in terms of its angular measure. Therefore no APD sequence exists.
\vskip2pt

\nid{\bf Proof.} From the lower bound in \eqref{bounds22A} we see that \eqref{Ahatlimit4} holds for the maximal configurations
$\{\^{\bf y}^q\}$. Because $\^{\bf y}^q$ minimizes the largest empty cap, \eqref{Ahatlimit4} holds for all nonzero sequences $\{\^{\bf y}^q\}$.\hfill$\square$
\vskip4pt

Lastly, the standardized limits of the univariate marginal distributions are as follows: for each fixed $i\ge1$, \begin{align}
\ts\sqrt{\frac{12}{q^2-1}}(\~{\bf U}_{\bar{\bf y}^q}^{q-2})_i&\tod\mathrm{Uniform}\big(-\sqrt{3},\sqrt{3}\,\big);\label{A}\\
\ts \sqrt{\frac{12}{q^2-1}}(\~{\bf U}_{\^{\bf y}^q}^{q-2})_i &\ts=O_p\big(\frac{t_2}{\sqrt{\log q}}\big)\top0;\label{B}\\
\ts\sqrt{\frac{12}{q^2-1}}(\~{\bf U}_{\u{\bf y}^q}^{q-2})_i&\tod N(0,1);\label{C}\\
\ts\sqrt{\frac{12}{q^2-1}}(\~{\bf U}_{\|\bar{\bf y}^q\|}^{q-2})_i&\tod N(0,1).\label{D}
\end{align}

\begin{table}[h!]
\centering
\begin{tabular}{ |p{2.6cm}||p{2.8cm}|p{2.8cm}|p{.8cm}|p{.8cm}|p{.8cm}|p{.8cm}|  }
\hline
  ${\bf y}^q$&$\lim\limits_{q\to\infty}L^{q-2}(\Pi({\bf y}^q))$&$\lim\limits_{q\to\infty}A^{q-2}(\Pi({\bf y}^q))$  &APF&APU&APE&$\!\!N\!(0,\!1)$\\
 \hline
 \hline
$\bar{\bf y}^q$  regular  &$1-\Phi(\sqrt{3})$&$\pi/2$&no&no&yes&no\\
 \hline
$\^{\bf y}^q$  maximal  &0&$\pi/2$&yes&no&no&no \\
 \hline
$\u{\bf y}^q$ normal &$1/2$&$\pi/2$&no&no&yes&yes    \\
\hline
$\|\bar{\bf y}^q\|$  spherical & ``0"&$``\pi/2"$&``yes"&``yes" &``no"&``yes"\\
 \hline
\end{tabular}
\caption{The first three rows refer to the discrete uniform distribution on the permutations in $\Pi({\bf y}^q)$. The fourth row refers to the continuous uniform distribution on the sphere $\~{\cal S}_{\|\bar{\bf y}^q\|}^{q-2}$, where the ``entries" hold trivially. The second and third columns  show the limiting LECDs and LECADs, respectively. The final column indicates whether or not the univariate marginal distributions converge to $N(0,1)$, a necessary condition for APU.}
\label{table:1}
\end{table}

Our asymptotic results for the LECDs, LECADs, and univariate marginal distributions of the regular, maximal, and normal configurations are summarized in Table 1. Neither the regular nor normal sequences is APU, nor is the maximal sequence APU even though it is APF. Therefore we conjecture, albeit somewhat weakly, that the answer to the following question is no:
\vskip2pt

\nid {\it Question 4: Does any  APU sequence $\{{\bf y}^q\in\~{\cal M}_\le^{q-1}\}$ exist?} 
\vskip10pt

\begin{table}[h!]
\centering
\begin{center}
\begin{tabular}{ |c|c|c|c|c|c|c|c|c|c|c|c| }
 \hline
 $q$&$\bar{\bf y}^q$&$\^{\bf y}^q$&$\u{\bf y}^q$\\
 \hline
 3&  (0,\,1)&(0,\,1)&(0,\,1)\\
4&(.5,\,1.5)&(.242,\,1.56)&(.459,\,1.51)\\
5&(0,\,1,\,2)&(0,\,.490,\,2.18)&(0,\,.909,\,2.04)\\
6&(.5,\,1.5,\,2.5)&(.219,\,.756,\,2.85)&(.436,\,1.37,\,2.59)\\
 \hline
\end{tabular}
\end{center}
\caption{The regular, maximal, and normal configurations for $q=3,4,5,6$. The $q$ components of each vector ${\bf y}^q$ are symmetric about 0 so only the nonnegative components are shown.}
\label{table:2}
\end{table}

Some exact values of $\bar{\bf y}^q$, $\^{\bf y}^q$, and $\u{\bf y}^q$, are shown in Table 2. 
For $q=3$, $\bar{\bf y}^3=\^{\bf y}^3=\u{\bf y}^3$, while for $q\ge4$ the components of $\^{\bf y}^q$ disperse more rapidly than those of $\bar{\bf y}^q$ and $\u{\bf y}^q$ as $q$ increases. This is also seen from the following asymptotic comparisons of the magnitudes of the ranges of the univariate marginal distributions: for each $i=1,\dots,q$,
\begin{align*}
\ts\big|\mathrm{range}[(\~{\bf U}_{\bar{\bf y}^q}^{q-2})_i]\big|&=\ts\big|[-\frac{q-1}{2},\,\frac{q-1}{2}]\big|&&\!\!\!=q&&\!\!\!=O(q)\\
\big|\mathrm{range}[(\~{\bf U}_{\^{\bf y}^q}^{q-2})_i]\big|&=\ts\big|\big[-\frac{\|\bar{\bf y}^q\|}{\|\^{\bf a}^q\|}\^a_q^q,\,\frac{\|\bar{\bf y}^q\|}{\|\^{\bf a}^q\|}\^a_q^q\big]\big|&&\!\!\!=\ts\frac{q-1}{\|\^{\bf a}^q\|}&&\!\!\!=\ts O\Big(\frac{q^{\frac{3}{2}}}{\sqrt{\log q}}\Big)\\
\big|\mathrm{range}[(\~{\bf U}_{\u{\bf y}^q}^{q-2})_i]\big|&=\ts\big|\big[-\frac{\|\bar{\bf y}^q\|}{\|\u{\bf a}^q\|}\u a_q^q,\,\frac{\|\bar{\bf y}^q\|}{\|\u{\bf a}^q\|}\u a_q^q\big]\big|&&\!\!\!\ts\sim\sqrt{\frac{q^2-1}{3}}\Phi^{-1}(\frac{q}{q+1})&&\!\!\!=O\Big(q\sqrt{\log q}\Big)\\
\big|\mathrm{range}[(\~{\bf U}_{\|\bar{\bf y}^q\|}^{q-2})_i]\big|&=\ts\big|[-\|\bar{\bf y}^q\|,\,\|\bar{\bf y}^q\|]\big|&&\!\!\!=\ts\sqrt{\frac{q(q^2-1)}{3}}&&\!\!\!=O\big(q^{\frac{3}{2}}\big).
\end{align*}

The four ranges satisfy 
\begin{equation}\label{rangeorder}
\mathrm{regular} \ll \mathrm{normal} \ll \mathrm{maximal} \ll \mathrm{spherical},
\end{equation}
where ``$\ll$" indicates $o(\cdot)$, whereas the limiting distributions of the univariate marginals in \eqref{A}-\eqref{D} satisfy
\begin{equation}\label{marginalorder}
\mathrm{maximal} \ll_p \mathrm{regular} \approx_p \mathrm{normal} \approx_p \mathrm{spherical},
\end{equation}
where ``$\ll_p$" indicates $o_p(\cdot)$ and ``$\approx_p$" indicates $O_p(\cdot)$. The ordering \eqref{marginalorder} is somewhat unexpected since the maximal configuration is the only one of the three uniform permutation distributions that is APF.
\vskip6pt

\nid{\bf 7. The regular, maximal, and normal permutohedra.} 

\nid The {\it regular permutohedron}\footnote{\label{Foot12} a.k.a. permutahedron.} $\mathfrak{R}^q$ is defined to be the convex hull of $\Pi(\bar{\bf x}^q)$, the set of all $q!$ permutations of the regular configuration $\bar{\bf x}^q\equiv(1,2,\dots,q)'$. It is a convex polyhedron in ${\cal M}_{\bar{\bf x}^q}^{q-1}$ (cf. \eqref{Mtilde}) of affine dimension $q-1$. Equivalently we shall consider the congruent polyhedron $\~{\frak R}^q\equiv\Ome_q\,\mathfrak{R}^q$, the translation of $\mathfrak{R}^q$ into $\~{\cal M}^{q-1}$, so
$\~{\frak R}^q$ is the convex hull of  $\Pi(\bar{\bf y}^q)$ (cf. \eqref{yqdef}). Thus the uniform distribution $\~{\bf U}_{\bar{\bf y}^q}^{q-2}$ is the uniform distribution on the vertices of $\~{\frak R}^q$.

Proposition 3.2 shows that  $\Pi(\bar{\bf y}^q)$ occupies a vanishingly small portion of the sphere $\~{\cal S}_{\|\bar{\bf y}^q\|}^{q-2}$ as $q\to\infty$. Similarly, it will now be shown that $\~{\frak R}^q$ occupies a vanishingly small portion of the corresponding ball $\~{\frak B}^q:=\~{\frak B}_{\|\bar{\bf y}^q\|}^{q-2}$ in which $\~{\frak R}^q$ is inscribed. 


\vskip4pt

\nid{\bf Proposition 7.1.} As $q\to\infty$, $\frac{\mathrm{Vol}(\~{\frak R}^q)}{\mathrm{Vol}(\~{\frak B}^q)}\to0$ at a geometric rate.
\vskip2pt

\nid{\bf Proof.} From Proposition 2.11 of Baek and Adams [BA] with $d=q-1$, the volume of $\~{\frak R}^q$ is $q^{q-\frac{3}{2}}$, while the volume of $\~{\frak B}^q$ is
\begin{equation*}
\ts\frac{\pi^{\frac{q-1}{2}}\|{\bf y}^q\|^{q-1} }{\Gam(\frac{q+1}{2})}=\frac{\pi^{\frac{q-1}{2}} }{\Gam(\frac{q+1}{2})}
\big[\frac{q(q^2-1)}{12}\big]^{\frac{q-1}{2}}.
\end{equation*}
Therefore, using Stirling's formula, the ratio of the volumes is given by
\begin{align}
\frac{\mathrm{Vol}(\~{\frak R}^q)}{\mathrm{Vol}(\~{\frak B}^q)}&=\ts\big(\frac{12}{\pi}\big)^{\frac{q-1}{2}}\frac{q^{q-\frac{3}{2}}\Gam(\frac{q+1}{2})}{[q(q^2-1)]^{\frac{q-1}{2}}}\nonumber\\
  &\sim\ts\big(\frac{\pi}{e}\big)^{\frac{1}{2}} \big(\frac{6}{\pi e}\big)^{\frac{q-1}{2}}\nonumber \\
   &\approx\ts1.0750\,\big(0.7026\big)^{\frac{q-1}{2}}\hskip80pt\label{Volratio2} 
\end{align}
as $q\to\infty$, which converges to zero at a geometric rate.
\vskip4pt

\nid{\bf Remark 7.2.} By comparison, the cube $\~{\frak C}^q$ inscribed in $\~{\frak B}^q$ has vertices
\begin{equation*}
\ts\big(\pm\frac{\|y^q\|}{\sqrt{q}},\dots,\pm\frac{\|y^q\|}{\sqrt{q}}\big)=\Big(\pm\sqrt{\frac{q^2-1}{12}},\dots,\pm\sqrt{\frac{q^2-1}{12}}\,\Big),
\end{equation*}
so
\begin{align}
\frac{\mathrm{Vol}(\~{\frak C}^q)}{\mathrm{Vol}(\~{\frak R}^q)}&=\ts\big(\frac{q^2-1}{3}\big)^{\frac{q-1}{2}}q^{\frac{3}{2}-q}\nonumber\\
&\sim\ts\frac{q^{\frac{1}{2}}}{3^{\frac{q-1}{2}}}\nonumber\\
 &\approx\ts q^{\frac{1}{2}}\,(0.3333)^{\frac{q-1}{2}}\label{sim2}
\end{align}
as $q\to\infty$, which also converges to zero at a geometric rate. Therefore 
\begin{equation}\label{Volcompare}
{\mathrm{Vol}(\~{\frak C}^q)}\ll{\mathrm{Vol}(\~{\frak R}^q)}\ll{\mathrm{Vol}(\~{\frak B}^q)}
\end{equation}
for large $q$.\hfill$\square$
\vskip4pt

Next, define the {\it maximal permutohedron} $\~{\frak M}^q$ ({\it normal permutohedron} $\~{\frak N}^q$) to be the convex hull of $\Pi(\^{\bf y}^q)$ ($\Pi(\u{\bf y}^q))$, the set of all $q!$ permutations of the maximal configuration $\^{\bf y}^q$ (normal configuration $\u{\bf y}^q$). Like the regular permutohedron $\~{\frak R}^q$ defined in \S7, $\~{\frak M}^q$ and $\~{\frak N}^q$ are convex polyhedrons in $\~{\cal M}^{q-1}$ (cf. \eqref{Mtilde}) of affine dimension $q-1$. Thus the uniform distribution $\~{\bf U}_{\^{\bf y}^q}^{q-2}$ ($\~{\bf U}_{\u{\bf y}^q}^{q-2}$) is the uniform distribution on the vertices of $\~{\frak M}^q$ ($\~{\frak N}^q$). the following question is suggested:
\vskip4pt

\nid{\it Question 5: What are the volumes of $\~{\frak M}^q$ and $\~{\frak N}^q$? As in Proposition 7.1 and Remark 7.2, compare ${\mathrm{Vol}(\~{\frak R}^q)}$, ${\mathrm{Vol}(\~{\frak M}^q)}$, ${\mathrm{Vol}(\~{\frak N}^q)}$, and ${\mathrm{Vol}(\~{\frak B}^q)}$.}

\nid We conjecture, again somewhat weakly, that as $q\to\infty$, 
\begin{equation}\label{Volcompare1}
{\mathrm{Vol}(\~{\frak R}^q)}\ll{\mathrm{Vol}(\~{\frak M}^q)}\ll{\mathrm{Vol}(\~{\frak B}^q)},
\end{equation}
more precisely, that $\frac{\mathrm{Vol}(\~{\frak R}^q)}{\mathrm{Vol}(\~{\frak M}^q)}\to0$ at a geometric rate and $\frac{\mathrm{Vol}(\~{\frak M}^q)}{\mathrm{Vol}(\~{\frak B}^q)}\to0$ at a slower rate. Similar results are expected if $\~{\frak M}^q$ is replaced by $\~{\frak N}^q$.
\vskip6pt

\nid{\bf 8. Concluding remarks.}

\nid We conclude with a final question and remark.
\vskip4pt

\nid {\it Question 6: If the permutation group is replaced by some other finite subgroup $G$ of orthogonal transformations on $\mathbb{R}^q$, how close to spherical uniformity is the $G$-orbit $\Pi_G({\bf y}^q)\equiv\{g{\bf y}^q\mid g\in G\}$ for nonzero ${\bf y}^q\in\mathbb{R}^q$?} 

\nid Finite reflection groups (Coxeter groups) acting on $\mathbb{R}^q$ for all $q\ge2$ are of particular interest, cf. [EP], [GB]. These include, and in fact are limited to, the permutation (= symmetric) group, the alternating group, and the group generated by all permutations and sign changes of coordinates.
\vskip4pt

\nid{\bf Remark 8.1.}  In coding theory, a finite set of $N$ points on a $d$-sphere is called a {\it spherical code}, cf. Leopardi [L1], [L2]. A question of major interest is the construction of spherical codes having small spherical discrepancy for large $N$ with $d$ held fixed (recall Definition 2.2). Thus our sets $\Pi({\bf x})$ and $\Pi({\bf y})$, consisting of all q! permutations of ${\bf x}$ and ${\bf y}$, can be viewed as spherical codes of a special type; we suggest that these be called {\it permutation codes}. We are interested in a similar question: which if any permutation codes are APU, that is, have small spherical discrepancy as $q\to\infty$?  Here, however, $N=q!$ and $d=q-2$, so both $N\to \infty$ and $d\to\infty$ in our case.   \hfill$\square$
\vskip10pt


\nid{\bf Appendix. Subindependence of coordinate halfspaces.}

\nid  The following inequality was used in the proof of Proposition 3.2:
\vskip4pt

\nid{\bf Proposition A.1.} Let ${\bf U}_n\equiv(U_1,\dots,U_n)'$ be uniformly distributed on the unit $(n-1)$-sphere ${\cal S}^{n-1}$ in ${\bf R}^n$. For any positive real numbers $t_1,\dots,t_n$, 
\begin{equation}\label{subindep}
\Pr[\cap_{i=1}^n\{U_i\le t_i\}]\le\prod\nolimits_{i=1}^n\Pr[U_i\le t_i].
\end{equation}
\nid{\bf Proof.} The proof is modelled on that of Proposition 2.10 in Barthe {\it et al.} [BGLR]. We shall show more generally that for $1\le r<n$,
\begin{equation}\label{negassoc}
\Pr[\cap_{i=1}^n\{U_i\le t_i\}]\le\Pr[\cap_{i=1}^r\{U_i\le t_i\}]\Pr[\cap_{i=r+1}^n\{U_i\le t_i\}].
\end{equation}

Because ${\bf U}_n$ is the unique orthogonally invariant distribution on ${\cal S}^{n-1}$, 
\begin{equation*}
{\bf U}_n\equiv\begin{pmatrix}{\bf U}_r\\{\bf U}_{-r}\end{pmatrix}\eqd\begin{pmatrix}\Psi&0\\0&I_{n-r}\end{pmatrix}\begin{pmatrix}{\bf U}_r\\{\bf U}_{-r}\end{pmatrix}=\begin{pmatrix}\psi{\bf U}_r\\{\bf U}_{-r}\end{pmatrix}
\end{equation*}
for every orthogonal $r\times r$ matrix $\Psi$,
where
\begin{align*}
{\bf U}_r&=(U_1,\dots,U_r)',\\
{\bf U}_{-r}&=(U_{r+1},\dots,U_n)'.
\end{align*}
Therefore
\begin{equation*}
\begin{pmatrix}{\bf W_r}\\{\bf U}_{-r}\end{pmatrix}\eqd\begin{pmatrix}\psi{\bf W_r}\\{\bf U}_{-r}\end{pmatrix},
\end{equation*}
where
\begin{align*}
{\bf W}_r\equiv(W_1,\dots,W_r)'&=\ts\frac{{\bf U}_r}{\|{\bf U}_r\|}\in{\cal S}^{r-1},\\
{\bf W}_{-r}\equiv(W_{r+1},\dots,W_n)'&=\ts\frac{{\bf U}_{-r}}{\|{\bf U}_{-r}\|}\in{\cal S}^{n-r+1}.
\end{align*}
Thus the conditional distribution of ${\bf W}_r| {\bf U}_{-r}$ is the same as that of $\Psi{\bf W}_r|{\bf U}_{-r}$ so, by uniqueness, is the uniform distribution on ${\cal S}^{r-1}$. Therefore ${\bf W}_r$ is independent of ${\bf U}_{-r}$, hence ${\bf W}_r$ is independent of $({\bf W}_{-r},\|{\bf U}_{-r}\|)$. Similarly, ${\bf W}_{-r}$ is independent of  $({\bf W}_r,\|{\bf U}_r\|)$. However, $\|{\bf U}_r\|$ and $\|{\bf U}_{-r}\|$ are statistically equivalent because $\|{\bf U}_r\|^2+\|{\bf U}_{-r}\|^2=\|{\bf U}_n\|^2=1$, hence ${\bf W}_{-r}$ is independent of $\|{\bf U}_{-r}\|$, so ${\bf W}_r$, ${\bf W}_{-r}$, and $\|{\bf U}_{-r}\|$ are mutually independent. Thus ${\bf W}_r$, ${\bf W}_{-r}$, and $\|{\bf U}_r\|$ are mutually independent, so
\begin{align*}
\Pr[\cap_{i=1}^n\{U_i&\le t_i\}]\\
=&\;\E\Big\{\Pr[\cap_{i=1}^r\{W_i\le t_i\|{\bf U}_r\|^{-1}\}\mid \|{\bf U}_r\|]\\
 &\ \ \cdot\Pr[\cap_{i=r+1}^n\{W_i\le t_i(1-\|{\bf U}_r\|^2)^{-1/2}\}\mid \|{\bf U}_r\|]\Big\}\\
\le&\;\E\big\{\Pr[\cap_{i=1}^r\{W_i\le t_i\|{\bf U}_r\|^{-1}\}\mid \|{\bf U}_r\|]\big\}\\
&\ \ \cdot\E\big\{\Pr[\cap_{i=r+1}^n\{W_i\le t_i(1-\|{\bf U}_r\|^2)^{-1/2}\}\mid \|{\bf U}_r\|]\big\}\\
=&\;\Pr[\cap_{i=1}^r\{U_i\le t_i\}]\cdot\Pr[\cap_{i=r+1}^n\{U_i\le t_i\}].
\end{align*}
The inequality holds because
\begin{equation*}
\Pr[\cap_{i=1}^r\{W_i\le t_i\|{\bf U}_r\|^{-1}\}\mid\|{\bf U}_r\|]
\end{equation*}
is decreasing in $\|{\bf U}_r\|$ while
\begin{equation*}
\Pr[\cap_{i=r+1}^n\{U_i\le t_i(1-\|{\bf U}_r\|^2)^{-1/2}\}\mid\|{\bf U}_r\|]
\end{equation*}
is increasing in $\|{\bf U}_r\|$.\hfill$\square$
\vskip4pt

\nid{\bf Remark A.2.} The inequality \eqref{subindep} is a one-sided version for coordinate halfspaces of a  two-sided inequality for symmetric coordinate slabs, where $|U_i|$ appears in place of $U_i$; see [BGLR] pp. 329-330 and the references cited therein. As in [BGLR], it is straightforward to extend Proposition A.1 to distributions on the unit sphere in $\ell_p$ for $1\le p<\infty$.\hfill$\square$
\bigskip

 \nid{\it Acknowledgement.} This paper was prepared with invaluable assistance from Steve Gillispie. Warm thanks are also due to Persi Diaconis, Art Owen, and Jens Praestgaard for very helpful discussions about random permutations, spherical geometry, and uniform distributions on groups.
\bigskip
\bigskip

\centerline{\bf References}
\vskip8pt


\nid [AZ] Alishahi, K. and Zamani, M. S. (2015). The spherical ensemble and uniform distribution of points on the sphere. {\it Electronic J. Prob.} {\bf 20} 1-27.

\nid [BA] Baek, J. and Adams, A. (2009).  Some useful properties of the permutohedral lattice for Gaussian filtering. http://graphics.stanford.edu/papers /permutohedral/permutohedral\_techreport.pdf.

\nid [BP] Ball, K. and Perissinaki, I. (1998). The subindependence of coordinate slabs in $\ell_p^n$ balls. {\it Israel J. Math.} {\bf 107} 289-299.

\nid [BGLR] Barthe, F., Gamboa, F., Lozada-Chang, L., Rouault, A. (2010). Generalized Dirichlet distributions on the ball and moments. {\it Latin American J. Prob. Math. Statist.} {\bf 7} 319-340.


\nid [DEOPSS] Das Gupta, S., Eaton, M. L., Olkin, I., Perlman, M. D., Savage, L. J., Sobel, M. (1972). Inequalities for the probability content of convex regions for elliptically contoured distributions. {\it Proc. Sixth Berkely Symp. Math. Statist. Prob.} {\bf 2} 241-265.




\nid [E] Eaton, M. L. (1989). {\it Group Invariance Applications in Statistics.} Regional Conference Series in Probability and Statistics Vol. 1, Institute of Mathematical Statistics.


\nid [EP] Eaton, M. L. and Perlman, M. D. (1977). Reflection groups, generalized Schur functions, and the geometry of majorization. {\it Ann. Prob.} {\bf 5} 829-860.

\nid [FS] Fung, T. and Seneta, E. (2018). Quantile function expansion using regularly varying functions. {\it Method. Comp. Appl. Probab.} {\bf 20} 1091-1103.

\nid [GB] Grove, L. C. and Benson, C. T. (1985). {\it Finite Reflection Groups}, 2nd Ed. Springer, New York.





\nid [L1] Leopardi, P. (2007). Distributing points on the sphere: partitions, separation, quadrature and energy. Ph.D.
thesis, The University of New South Wales (2007).

\nid [L2] Leopardi, P. (2013). Discrepancy, separation and Riesz energy of finite point sets on the unit sphere. {\it Adv. Comp. Math.} {\bf 39} 27-43.


\nid [MO] Marshall, A. W. and Olkin, I. (1979). {\it Inequalities: Theory of Majorization and Its Applications.} Academic Press, New York.

\nid [QG] Qi, F. and Guo, B.-N. (2011). Sharp bounds for harmonic numbers. {\it Applied Math. and Computation} {\bf 218} 991–995.




\nid [W] Wendel, J. G. (1948). Note on the gamma function. {\it Amer. Math. Monthly} {\bf 55} 563-564.



\bigskip

\newpage

\end{document}